\documentclass[12pt,reqno]{amsart}

\usepackage{geometry}
\usepackage{amsmath,amssymb,mathtools}
\usepackage{mathrsfs,eucal}
\usepackage{bm}

\usepackage{graphicx}
\usepackage{multirow,bigdelim}

\usepackage{enumitem}

\usepackage[colorlinks=true, citecolor=blue, linkcolor=blue, urlcolor=blue]{hyperref}

\allowdisplaybreaks[1]

\makeatletter
\@namedef{subjclassname@2020}{\textup{2020} Mathematics Subject Classification}
\makeatother

\providecommand{\doi}[1]{}
\DeclareUrlCommand\doifont{\urlstyle{tt}}
\renewcommand{\doi}[1]{%
  \href{https://doi.org/#1}{\doifont{https://doi.org/#1}}}

\newtheorem{lemma}{Lemma}[section]
\newtheorem{theorem}[lemma]{Theorem}
\newtheorem{proposition}[lemma]{Proposition}
\newtheorem{corollary}[lemma]{Corollary}
\newtheorem{question}[lemma]{Question}
\theoremstyle{definition}
\newtheorem{definition}[lemma]{Definition}
\newtheorem{example}[lemma]{Example}
\newtheorem{remark}[lemma]{Remark}

\numberwithin{equation}{section}

\newcommand{\bdf}{\begin{definition}}
\newcommand{\edf}{\end{definition}}
\newcommand{\blem}{\begin{lemma}}
\newcommand{\elem}{\end{lemma}}
\newcommand{\bthm}{\begin{theorem}}
\newcommand{\ethm}{\end{theorem}}
\newcommand{\bpf}{\begin{proof}}
\newcommand{\epf}{\end{proof}}
\newcommand{\bprop}{\begin{proposition}}
\newcommand{\eprop}{\end{proposition}}
\newcommand{\bcor}{\begin{corollary}}
\newcommand{\ecor}{\end{corollary}}
\newcommand{\brem}{\begin{remark}}
\newcommand{\erem}{\end{remark}}
\newcommand{\bquest}{\begin{question}}
\newcommand{\equest}{\end{question}}
\newcommand{\bex}{\begin{example}}
\newcommand{\eex}{\end{example}}

\newcommand{\benu}{\begin{enumerate}\renewcommand{\labelenumi}{{\rm (\arabic{enumi})}}\renewcommand{\itemsep}{0pt}}
\newcommand{\eenu}{\end{enumerate}}

\newcommand{\R}{\mathbb{R}}
\newcommand{\C}{\mathbb{C}}

\newcommand{\bB}{\mathbb{B}}

\newcommand{\cH}{\mathcal{H}}

\newcommand{\cL}{\mathcal{L}}
\newcommand{\cM}{\mathcal{M}}
\newcommand{\cN}{\mathcal{N}}

\newcommand{\cP}{\mathcal{P}}

\newcommand{\cV}{\mathcal{V}}

\newcommand{\e}{\varepsilon}

\DeclareMathOperator{\id}{id}
\DeclareMathOperator{\tr}{tr}
\DeclareMathOperator{\Tr}{Tr}

\newcommand{\ip}[1]{\mathopen{\langle}#1\mathclose{\rangle}}

\begin{document}
\title{Geometric mean and Lebesgue-type decomposition of completely positive maps}
\author{Rui OKAYASU}
\address{Department of Mathematics Education, Osaka Kyoiku University, Kashiwara, Osaka 582-8582, JAPAN}
\email{rui@cc.osaka-kyoiku.ac.jp}
\date{}
\thanks{The author was partially supported by JSPS KAKENHI Grant Number JP26K06843.}

\begin{abstract}
We introduce the geometric mean and the parallel sum of completely positive (CP) maps between von Neumann algebras, based on the Pusz--Woronowicz theory of positive sesquilinear forms. 
We provide a concrete characterization via a block matrix positivity condition and establish their fundamental properties, including the AM--GM--HM inequality with respect to the CP order. 

In finite-dimensional settings, our construction is compatible with the Choi--Jamio\l{}kowski correspondence, under which the geometric mean of CP maps corresponds to the Kubo--Ando geometric mean of their Choi matrices. 
This yields a natural operator-theoretic framework for interpolating quantum channels.

As an application, we obtain index-type inequalities for conditional expectations in subfactor theory. 

Finally, we establish a Lebesgue-type decomposition of CP maps via a parallel sum construction, thereby providing a unified framework that simultaneously generalizes Ando's decomposition of bounded positive operators and Kosaki's decomposition of normal positive functionals on von Neumann algebras.
\end{abstract}

\maketitle

\section{Introduction}

Completely positive (CP) maps play a central role in operator algebras and quantum information theory, serving as the natural morphisms between noncommutative spaces. 
Despite their fundamental importance, a systematic theory of operator means and Lebesgue-type decompositions for CP maps has remained relatively underdeveloped, in contrast to the well-established theories for bounded positive operators and normal positive functionals.

The aim of this paper is to develop a unified geometric and order-theoretic framework for CP maps that extends both operator mean theory in the sense of Kubo--Ando \cite{ka} and Lebesgue-type decomposition theory in the sense of Ando \cite{and2} and Kosaki \cite{kos_1985}. 
Our approach is based on the Pusz--Woronowicz theory \cite{pw1,pw2} of positive sesquilinear forms and the spatial realization of CP maps via the Stinespring dilation \cite{st}.

The starting point of our work is the observation that the geometric mean of positive sesquilinear forms provides a natural mechanism for interpolating between noncommutative structures without relying on multiplicative operations. 
This idea, originally introduced by Pusz and Woronowicz, has proved highly effective in mathematical physics, particularly in the study of relative entropy and quantum information geometry, as developed by Uhlmann \cite{uhl_1977} and Kosaki \cite{k1,kos_1982,kos_1986}. 
It is also closely related to subsequent developments in operator mean theory and noncommutative integration, for instance in the work of Hiai and Kosaki \cite{hk2}. 
In this paper, we show that this formalism extends naturally to CP maps and gives rise to a rich order-theoretic and geometric structure.

Given two CP maps $\Phi$ and $\Psi$ between von Neumann algebras, we define their geometric mean $\Phi \# \Psi$ through the associated positive sesquilinear forms. 
We establish a concrete characterization of this mean via a block matrix positivity condition, which yields an operator-inequality formulation analogous to the Kubo--Ando theory. 
In particular, we introduce the arithmetic, geometric, and harmonic means of CP maps and prove the AM--GM--HM inequality
\[
\Phi \, ! \, \Psi \le_{\mathrm{cp}} \Phi \# \Psi \le_{\mathrm{cp}} \Phi \nabla \Psi.
\]
We further show that these means satisfy fundamental structural properties such as joint monotonicity, concavity, and the transformer inequality.

In finite-dimensional settings, our construction is compatible with the Choi--Jamio\l kowski correspondence \cite{choi, jam}. Namely, the geometric mean of CP maps corresponds to the Kubo--Ando geometric mean of their Choi matrices. 
This provides a natural operator-theoretic framework for interpolating quantum channels and reveals a structural mechanism for extracting their common features. 
Related aspects of operator means for CP maps have recently been investigated by Frenkel, Mosonyi, Vrana, and Weiner \cite{FMVW} in the context of quantum hypothesis testing.

As an application to operator algebras, we study geometric means of conditional expectations. 
While the geometric mean of conditional expectations is not itself a conditional expectation in general, we show that it admits a natural bimodule structure over the intersection algebra. 
We also obtain index-type inequalities that connect our construction to subfactor theory \cite{izu, pp}.

Beyond the theory of means, we develop a theory of connections for CP maps that parallels the Kubo--Ando theory for bounded positive operators. 
This shows that many fundamental order-theoretic properties extend to the noncommutative setting of CP maps.

A central result of this paper is a canonical Lebesgue-type decomposition for CP maps, constructed via a parallel sum limit. 
This provides a unified framework that simultaneously generalizes Ando's decomposition of bounded positive operators and Kosaki's decomposition of normal positive functionals on von Neumann algebras. 

Our approach is based on Arveson's Radon--Nikodym theorem \cite{arv}, which allows us to realize CP maps in a common Stinespring space and to interpret absolute continuity and singularity in terms of support projections in the commutant. 

A key feature of our formulation is that it reduces analytic conditions---such as range inclusion, density, and closability---to simple geometric relations between closed subspaces. 
In particular, the seemingly different criteria appearing in the works of Ando and Kosaki are shown to be manifestations of a single underlying projection-theoretic structure.

Taken together, these results reveal a rich geometric and order structure on the space of CP maps, providing a unified perspective that connects operator means, noncommutative integration, and quantum information theory.

\section{Sesquilinear forms}

This section recalls the necessary parts of the Pusz--Woronowicz theory 
of positive sesquilinear forms \cite{pw1}, 
fixing notation for later use.

\subsection{Pusz--Woronowicz theory}

Let $\cV$ be a complex vector space.
We denote by $F(\cV)$ 
the set of all sesquilinear forms on $\cV$,
i.e., the set of all functions $\cV\times\cV\to\C$,
which are linear in the first and conjugate-linear
in the second variable.
Note that $F(\cV)$ naturally forms a complex vector space.
A sesquilinear form $\alpha\in F(\cV)$ is called {\em positive}
if 
$\alpha(x, x)\geq 0$
for all $x\in\cV$.
We denote by $F_+(\cV)$
the set of all positive sesquilinear forms on $\cV$.
Note that $F_+(\cV)$ is a convex cone in $F(\cV)$.
For $\alpha, \beta\in F(\cV)$,
we define a partial order $\alpha\geq \beta$ 
by  $\alpha-\beta\in F_+(\cV)$.

Let $\cH$ be a complex Hilbert space.
We denote by $\bB(\cH)$ the C$^*$-algebra of all bounded linear operators on $\cH$.
Let $h$ be a linear map from $\cV$ onto a dense subspace of $\cH$,
and let $A\in\bB(\cH)$.
We say that $\alpha\in F(\cV)$ is {\em represented by} $(h, A)$
if 
\[
\alpha(x, y)=\ip{Ah(x), h(y)}
\] 
for $x, y\in\cV$.
In this case, we write $\alpha\sim(h, A)$.
Note that $\alpha$ is positive if and only if
$A$ is positive.

If $\alpha\sim(h,A)$ and $\beta\sim(h,B)$ with $AB=BA$,
then we say that these representations are {\em compatible}. 
By \cite[Theorem 1.1]{pw1}, any pair of positive sesquilinear forms admits a compatible representation.

Let $f$ be a locally bounded Borel function on 
$\R_+^2=\{(r,s)\in\R^2 \mid r\geq 0, s\geq 0\}$. 
Recall that $f$ is {\em homogeneous}
if 
\[
f(\lambda r, \lambda s)=\lambda f(r, s)
\]
for $r, s, \lambda\in\R_+$.
The following theorem establishes the Pusz--Woronowicz functional calculus 
for positive sesquilinear forms.

\bthm[{\cite[Theorem 1.2]{pw1}}]\label{thm:pw1}
Let $f$
be a homogeneous locally bounded Borel function 
on $\R_+^2$.
For $\alpha, \beta\in F_+(\cV)$,
the sesquilinear form
\[
\gamma(x, y)=\ip{f(A,B)h(x), h(y)}
\]
is independent of the choice of compatible representations $\alpha\sim(h,A)$ and $\beta\sim(h, B)$.
\ethm

We denote by $f(\alpha, \beta)$ the resulting sesquilinear form $\gamma$ in Theorem \ref{thm:pw1}.
Note that the functions 
\[
f_G(r,s)=\sqrt{rs}
\quad\text{and}\quad
f_P(r,s)=\frac{rs}{r+s}\ (\text{with the convention $0/0=0$})
\]
are such homogeneous functions.
Therefore, for $\alpha, \beta\in F_+(\cV)$,
we define the {\em geometric mean} $\alpha\#\beta$
and
the {\em parallel sum} $\alpha : \beta$ by 
\[
\alpha\#\beta\coloneqq f_G(\alpha, \beta)
\quad\text{and}\quad
\alpha : \beta\coloneqq f_P(\alpha, \beta),
\]
respectively.
We also define the {\em harmonic mean} $\alpha\, !\, \beta$ by
\[
\alpha\, !\, \beta\coloneqq 2(\alpha : \beta).
\]

Let $\alpha, \beta, \gamma\in F_+(\cV)$.
We say that $\gamma$ is {\em dominated by}
$\{\alpha, \beta\}$ if
\[
|\gamma(x,y)|^2\leq\alpha(x,x)\beta(y,y)
\]
for $x,y\in\cV$.

Finally, we conclude this section with the following characterizations of the geometric mean and the parallel sum:

\bthm[{\cite[Theorem 2.1]{pw1}}]\label{thm:pw2}
Let $\alpha, \beta\in F_+(\cV)$.
Then
\[
\alpha\#\beta
=
\max\{
\gamma\in F_+(\cV)
\mid
\text{$\gamma$ is dominated by
$\{\alpha, \beta\}$}
\}.
\]
\ethm

\blem[{\cite[Lemma]{pw1}}]\label{lem:pw}
Let $\alpha, \beta\in F_+(\cV)$.
Then
\[
(\alpha : \beta)(z, z)=
\inf\{
\alpha(x,x)+\beta(y,y)\mid z=x+y, x,y\in\cV
\}
\quad\text{for}\ z\in\cV.
\]
\elem

\subsection{Matrices of sesquilinear forms}

Given $\gamma_{ij}\in F(\mathcal V)$
for $1\leq i,j\leq n$, define the matrix form 
$\Gamma=[\gamma_{ij}]$ on $\mathcal V^n$ by
\[
\Gamma(x,y)=\sum_{i,j}\gamma_{ij}(x_j,y_i)
\]
for $x=(x_1,\dots,x_n), y=(y_1,\dots,y_n)\in\cV^n$.
Note that $\Gamma(x, x)\geq 0$
if and only if 
the scalar matrix
$[\gamma_{ij}(x_j,x_i)]_{1\leq i,j\leq n}$ is 
positive.
Then we can obtain the block matrix characterization of the geometric mean of positive sesquilinear forms.

\bprop\label{prop:geometric}
Let $\alpha, \beta\in F_+(\cV)$.
Then
\[
\alpha\# \beta
=
\max
\left\{
\gamma\in F_+(\cV)
\;\middle|\;
\begin{bmatrix}
\alpha & \gamma \\
\gamma & \beta \\
\end{bmatrix}\geq 0
\right\}.
\]
\eprop

\bpf
Let $\gamma\in F_+(\cV)$.
The  proposition follows from Theorem \ref{thm:pw2} and the following fact:
\[
 \begin{bmatrix}
\alpha & \gamma \\
\gamma & \beta \\
\end{bmatrix}\geq 0
\quad
\text{if and only if}
\quad
|\gamma(x,y)|^2\leq\alpha(x, x)\beta(y, y)
\quad\text{for}\ x,y\in\cV.
\]
\epf

Similarly, one can derive block matrix characterizations 
of the parallel sum and the harmonic mean 
from Lemma~\ref{lem:pw}, 
using the variational formula for the parallel sum 
and a polarization argument.

\bprop\label{prop:harmonic}
Let $\alpha, \beta\in F_+(\cV)$. 
Then
\benu
\item
${\displaystyle
\alpha : \beta=\max
\left\{
\gamma\in F_+(\cV)  \;\middle|\;
\begin{bmatrix}
\alpha & \alpha \\
\alpha & \alpha+\beta \\
\end{bmatrix}
\geq
\begin{bmatrix}
\gamma & 0 \\
0 & 0 \\
\end{bmatrix}
\right\},
}$\vspace{2mm}
\item 
${\displaystyle
\alpha\, !\, \beta=\max
\left\{
\gamma\in F_+(\cV) \;\middle|\;
\begin{bmatrix}
2\alpha & 0 \\
0& 2\beta \\
\end{bmatrix}
\geq
\begin{bmatrix}
\gamma & \gamma \\
\gamma & \gamma \\
\end{bmatrix}
\right\}.}$
\eenu
\eprop

\section{Geometric mean and Parallel sum of completely positive maps}

Let $\cM$ and $\cN$ be von Neumann algebras. To apply the Pusz--Woronowicz theory to completely positive (CP) maps from $\cM$ to $\cN$, we use the standard form of the target algebra $\cN$.

We review the theory of bimodules over von Neumann algebras. For further details and background, see \cite{oot, pop, takebook}.

Let $L^2(\cN)$ be the standard form
of $\cN$, and $J$ be the modular conjugation. 
Note that $L^2(\cN)$ has the $\cN$-$\cN$ bimodule structure with left and right actions given by 
\[
x\xi y\coloneqq xJ y^* J\xi
\]
for $x, y\in\cN$ and $\xi\in L^2(\cN)$.
Let us consider the $n\times n$ matrix algebra $M_n(\C)$
with the canonical normalized trace $\tr$.
If we define the inner product on $M_n(\C)$ by
\[
\ip{x, y}\coloneqq\tr(y^*x)
\]
for $x, y\in M_n(\C)$,
then $M_n(\C)$ can be regarded as a Hilbert space.
Moreover, $M_n(\C)$ acts on itself from the left
such that the modular conjugation is the canonical involution
$J_{\tr}\colon x\to x^*$.
We also denote by $J^{(n)}$
the modular conjugation on
the standard form of $M_n(\cN)=\cN\otimes M_n(\C)$.

We denote by $\mathrm{CP}(\cM, \cN)$ 
the set of all CP maps from $\cM$ to $\cN$.
For $\Phi,\Psi\in\mathrm{CP}(\cM, \cN)$,
we define a partial order
$\Phi\geq_{\mathrm{cp}}\Psi$
if $\Phi-\Psi$ is completely positive.

For $\Phi\in\mathrm{CP}(\cM, \cN)$, 
we define a positive sesquilinear form $s_\Phi$
on the algebraic tensor product $\cM\odot L^2(\cN)$ by
\[
s_\Phi(x \otimes \xi, y \otimes \eta) \coloneqq \langle \Phi(y^* x) \xi, \eta \rangle
\]
for $x, y\in\cM$ and $\xi, \eta\in L^2(\cN)$.
Note that, 
for $\Phi, \Psi\in\mathrm{CP}(\cM, \cN)$,
$\Phi\leq_{\mathrm{cp}}\Psi$
if and only if $s_\Phi\leq s_\Psi$.
Via separation and completion,
we construct the Hilbert space
$\cH_\Phi$.
The left action of $\cM$ and the right action of $\cN$ are defined by
\[
a(x\otimes_\Phi \xi)b\coloneqq
(ax)\otimes_\Phi (\xi b)
\]
for $a, x\in \cM$, $b\in\cN$ and $\xi\in L^2(\cN)$.

We now give a short proof of \cite[Theorem 1.1]{pw1} in our setting.
For $\Phi, \Psi\in\mathrm{CP}(\cM, \cN)$,
we set 
$\alpha=s_\Phi$ and $\beta=s_\Psi$.

Put $\gamma=\alpha+\beta$.
For $b\in \cN$, we have
\begin{equation}\label{eq:right}
\gamma(x\otimes (\xi b), y\otimes \eta)
=
\gamma(x\otimes \xi, y\otimes (\eta b^*)).
\end{equation}
Via separation and completion of $\cM\odot L^2(\cN)$,
we construct a Hilbert space $\cH_\gamma$.
Denote by $h$ the canonical map from $\cM \odot L^2(\cN)$ to $\cH_\gamma$.
Then there exist bounded positive operators $A, B$ on $\cH_\gamma$ such
that
\[
\alpha(x\otimes \xi, y\otimes \eta)
=
\ip{
Ah(x\otimes \xi), h(y\otimes \eta)
}_\gamma
\]
and
\[
\beta(x\otimes \xi, y\otimes \eta)
=
\ip{
Bh(x\otimes \xi), h(y\otimes \eta)
}_\gamma.
\]
Note that $A+B=I$.
In particular, $A$ and $B$ commute.
For $b\in \cN$,
we can define a bounded operator $\rho(b)$ on $\cH_\gamma$ by
\[
\rho(b)h(x\otimes \xi)=h(x\otimes (\xi b)),
\]
for $x\in\cM$ and $\xi\in L^2(\cN)$,
because
\begin{align*}
\ip{h(x\otimes (\xi b)), h(x\otimes (\xi b))}_\gamma
&= \alpha(x\otimes (\xi b), x\otimes (\xi b)) + \beta(x\otimes (\xi b), x\otimes (\xi b)) \\
&\leq \|b\|^2\Big(\alpha(x\otimes \xi, x\otimes \xi) + \beta(x\otimes \xi, x\otimes \xi)\Big) \\
&= \|b\|^2\ip{h(x\otimes \xi), h(x\otimes \xi)}_\gamma.
\end{align*}
By using (\ref{eq:right}), we obtain
$\rho(b)^*=\rho(b^*)$. 
Moreover,
\begin{align*}
\ip{A\rho(b)h(x\otimes \xi), h(y\otimes \eta)}_\gamma
&= \ip{Ah(x\otimes (\xi b)), h(y\otimes \eta)}_\gamma \\
&= \alpha(x\otimes (\xi b), y\otimes \eta) \\
&= \alpha(x\otimes \xi, y\otimes (\eta b^*)) \\
&= \ip{Ah(x\otimes \xi), \rho(b^*)h(y\otimes \eta)}_\gamma \\
&= \ip{\rho(b)Ah(x\otimes \xi), h(y\otimes \eta)}_\gamma.
\end{align*}
Thus 
$\rho(b)A=A\rho(b)$. 
Similarly, we obtain $\rho(b)B=B\rho(b)$. 
Hence, 
for any homogeneous locally bounded Borel function $f$
on $\R_+^2$,
we have
\[
f(A, B)\rho(b)=\rho(b)f(A, B).
\] 
Therefore, for $b\in\cN$,
\begin{equation}\label{eq:right-action}
f(s_\Phi, s_\Psi)(x\otimes (\xi b), y\otimes \eta)=
f(s_\Phi, s_\Psi)(x\otimes \xi, y\otimes (\eta b^*)).
\end{equation}
Similarly, we also have, for $a\in\cM$,
\begin{equation}\label{eq:left-action}
f(s_\Phi, s_\Psi)((ax)\otimes \xi, y\otimes \eta)=
f(s_\Phi, s_\Psi)(x\otimes \xi, (a^*y)\otimes \eta).
\end{equation}
These equations play a crucial role in recovering a CP map from the resulting sesquilinear form.

\subsection{Geometric mean}

Let $\Phi, \Psi\in\mathrm{CP}(\cM, \cN)$.
By Pusz--Woronowicz theory, we obtain the geometric mean $\gamma=s_\Phi\# s_\Psi$ on $\cM\odot L^2(\cN)$.
Since the positive sesquilinear form $\gamma$ is dominated by 
$\{s_\Phi, s_\Psi\}$,
for $x, y\in \cM$ and $\xi, \eta\in L^2(\cN)$,
we have
\begin{align*}
|\gamma(x\otimes \xi, y\otimes \eta)|^2
&\leq
s_\Phi(x\otimes \xi, x\otimes \xi)s_\Psi(y\otimes \eta, y\otimes \eta)
\\
&=
\ip{
\Phi(x^*x)\xi, \xi
}
\ip{
\Psi(y^*y)\eta, \eta
}.
\end{align*}
We define a Hilbert space $\cH_\gamma$ via separation and completion
of $\cM\odot L^2(\cN)$ with respect to $\gamma$.
We denote by $h\colon \cM\odot L^2(\cN)\to\cH_\gamma$ the canonical map.

We first claim that the Hilbert space $\cH_\gamma$ carries a left action of $\cM$ and a right action of $\cN$.
For a fixed $a\in\cM$, we define the positive sesquilinear form
$\gamma^a$ on $\cM\odot L^2(\cN)$ by
\[
\gamma^a(x\otimes \xi, y\otimes \eta)
\coloneqq
\gamma((ax)\otimes \xi, (ay)\otimes \eta)
\]
for $x, y\in\cM$ and $\xi, \eta\in L^2(\cN)$.
Since
\begin{align*}
|\gamma^a(x\otimes \xi, y\otimes \eta)|^2
&\leq
\ip{
\Phi(x^*a^*ax)\xi, \xi
}
\ip{
\Psi(y^*a^*ay)\eta, \eta
}
\\
&\leq
\|a\|^4s_\Phi(x\otimes \xi, x\otimes \xi)
s_\Psi(y\otimes \eta, y\otimes \eta),
\end{align*}
the positive sesquilinear form $\gamma^a$ is dominated by
$\{\|a\|^2s_\Phi, \|a\|^2s_\Psi\}$.
Hence,
\[
\gamma^a\leq(\|a\|^2s_\Phi)\#(\|a\|^2s_\Psi)
=\|a\|^2(s_\Phi\# s_\Psi)=\|a\|^2\gamma.
\]
Therefore, we can define the left action of $\cM$ on $\cH_\gamma$ by
\[
a(x\otimes \xi)\coloneqq (ax)\otimes \xi
\]
for $a, x\in\cM$ and $\xi\in L^2(\cN)$.
Similarly, we can also define the right action of $\cN$ on $\cH_\gamma$ by
\[
(x\otimes \xi)b\coloneqq x\otimes (\xi b)
\]
for $x\in\cM$, $b\in \cN$ and $\xi\in L^2(\cN)$.
By (\ref{eq:right-action}), 
for $b\in\cN$, we have
\[
\gamma(x\otimes (\xi b), y\otimes \eta)=
\gamma(x\otimes \xi, y\otimes (\eta b^*)).
\]

Next we claim that an operator
\[
V\colon L^2(\cN)\to\cH_\gamma, 
\quad
\xi\mapsto h(1\otimes \xi)
\]
is bounded.
Indeed,
\begin{align*}
|\gamma(1\otimes \xi, 1\otimes \xi)|^2
&\leq
\ip{
\Phi(1)\xi, \xi
}
\ip{
\Psi(1)\xi, \xi
}
\\
&\leq
\|\Phi\|\|\Psi\|\|\xi\|^4.
\end{align*}
Hence, we obtain the bounded operator $V$.
Then we define the CP map $\Theta$ on $\cM$ by
\[
\Theta(x)=V^*xV
\]
for $x\in\cM$.

We finally show that $\Theta(x)\in\cN$ for $x\in \cM$.
Since, for $b\in\cN$,
\begin{align*}
\ip{
\Theta(x)Jb^*J\xi, \eta
}
&=
\ip{
xV(\xi b), V\eta 
}
\\
&=
\gamma(x\otimes (\xi b), 1\otimes \eta)
\\
&=
\gamma(x\otimes \xi, 1\otimes (\eta b^*))
\\
&=
\ip{
xV\xi, V(\eta b^*) 
}
\\
&=
\ip{
\Theta(x)\xi, \eta b^*
}
\\
&=
\ip{
Jb^*J\Theta(x)\xi, \eta
},
\end{align*}
we obtain $\Theta(x)\in(\cN')'=\cN$.
By (\ref{eq:left-action}),
we also have
\[
\ip{
\Theta(y^*x)\xi, \eta
}
=
(s_\Phi \# s_\Psi)(x\otimes \xi, y\otimes \eta).
\]

\bdf
For $\Phi, \Psi\in\mathrm{CP}(\cM, \cN)$,
the {\em geometric mean} $\Phi\#\Psi$
is defined by
\[
\Phi\#\Psi\coloneqq \Theta\in\mathrm{CP}(\cM, \cN).
\]
\edf

By the definition,
one can also see that if $\Phi$ and $\Psi$ are normal,
then so is $\Phi\#\Psi$.
The following can be easily verified.

\bprop
Let $\Phi, \Psi\in\mathrm{CP}(\cM, \cN)$.
Then the following conditions hold:
\benu
\item $\Phi\#\Psi=\Psi\#\Phi$,
\item $\Phi\#\Phi=\Phi$,
\item $\|\Phi\#\Psi\|\leq\|\Phi\|^{1/2}\|\Psi\|^{1/2}$.
\eenu
\eprop

\brem
Our construction of the geometric mean for CP maps avoids the approximation procedures typically required in the Kubo--Ando framework such as considering $\lim_{\e \downarrow 0}(A + \e I)$ for non-invertible operators. By defining the geometric mean through the Pusz--Woronowicz theory of sesquilinear forms on the standard representation, we can directly handle general CP maps. 
\erem

\brem
Let $\cM=\C$ and $\cN=\bB(\cH)$,
where $\cH$ is an infinite-dimensional Hilbert space.
For each $A\in\bB(\cH)^+$,
we define a CP map $\Phi_A\colon \C\to\bB(\cH)$
by
\[
\Phi_A(z)\coloneqq z A
\]
for $z\in\C$.
Then the corresponding sesquilinear form
$s_{\Phi_A}$ on $\C\otimes\cH\cong\cH$
is given by
\[
s_{\Phi_A}(z\otimes\xi, w\otimes\eta)
=z\bar w\ip{
A\xi, \eta
}.
\]
For any $A, B\in\bB(\cH)^+$,
the geometric mean of the sesquilinear forms $s_{\Phi_A}$ and $s_{\Phi_B}$
satisfies
\[
s_{\Phi_A}\#s_{\Phi_B}(z\otimes\xi, w\otimes\eta)
=
z\bar w
\ip{
(A\# B)\xi, \eta
},
\]
where $A\# B$ denotes the geometric mean of bounded positive operators.
Consequently,
we obtain 
\[
\Phi_A\#\Phi_B=\Phi_{A\# B}.
\]
This shows that our definition of the geometric mean for CP maps is a natural extension of the classical operator version.
\erem

Next, we give a block matrix characterization of the geometric mean $\Phi\#\Psi$.
Let $\Phi, \Psi, \Theta\in\mathrm{CP}(\cM,\cN)$. 
We define the linear maps
$\Pi\colon\cM\to M_2(\cN)$
and
$\tilde{\Pi}\colon M_2(\cM)\to M_2(\cN)$
by
\[
\Pi(x)\coloneqq
\begin{bmatrix}
\Phi(x) & \Theta(x) \\
\Theta(x) & \Psi(x) \\
\end{bmatrix}\quad
\text{and}\quad
\tilde\Pi
\begin{bmatrix}
x_{11} & x_{12} \\
x_{21} & x_{22} \\
\end{bmatrix}\coloneqq
\begin{bmatrix}
\Phi(x_{11}) & \Theta(x_{12}) \\
\Theta(x_{21}) & \Psi(x_{22}) \\
\end{bmatrix},
\]
respectively. 
By {\cite[Remark 1.2]{haa}},
$\Pi$ is CP if and only if $\tilde\Pi$ is CP.

\bprop\label{prop:geo}
Let $\Phi, \Psi \in \mathrm{CP}(\cM, \cN)$. 
Then
\[
\Phi \# \Psi = \max
\left\{ \Theta \in \mathrm{CP}(\cM, \cN)
\,\middle|\,
\begin{bmatrix}
\Phi & \Theta \\
\Theta & \Psi \\
\end{bmatrix} \geq_{\mathrm{cp}} 0
\right\}.
\]
\eprop

\bpf
Suppose that $\Pi\geq_{\mathrm{cp}}0$.
Then $\tilde\Pi\geq_{\mathrm{cp}}0$.
Recall that $\tilde\Pi$ is CP if and only if $s_{\tilde\Pi}$ is positive.
Since $s_{\tilde\Pi}$ is positive,
for $\zeta_1, \zeta_2\in \cM\odot L^2(\cN)$, we have 
\[
|s_\Theta(\zeta_1, \zeta_2)|^2 \leq s_\Phi(\zeta_1, \zeta_1) s_\Psi(\zeta_2, \zeta_2).
\]
In fact, for simple tensors $\zeta_1=x_1\otimes \xi_1$
and $\zeta_2=x_2\otimes \xi_2$,
if
\[
\boldsymbol{x}=
\begin{bmatrix}
x_1 & x_2 \\
0 & 0 \\
\end{bmatrix}\quad\text{and}\quad
\boldsymbol{\xi} =
\begin{bmatrix}
\xi_1 & 0 \\
\xi_2 & 0 \\
\end{bmatrix},
\]
then
\begin{align*}
0\leq
s_{\tilde\Pi}(\boldsymbol{x} \otimes \boldsymbol{\xi}, \boldsymbol{x} \otimes \boldsymbol{\xi})
&=
\ip
{
\tilde\Pi(\boldsymbol{x}^*\boldsymbol{x})\boldsymbol{\xi},
\boldsymbol{\xi}
}
\\
&=
\ip{
\begin{bmatrix}
\Phi(x_1^*x_1) & \Theta(x_1^*x_2) \\
\Theta(x_2^*x_1) & \Psi(x_2^*x_2) \\
\end{bmatrix}
\begin{bmatrix}
\xi_1 & 0 \\
\xi_2 & 0 \\
\end{bmatrix},
\begin{bmatrix}
\xi_1 & 0 \\
\xi_2 & 0 \\
\end{bmatrix}
}
\\
&=
\ip{
\Phi(x_1^*x_1)\xi_1, \xi_1
}
+ \ip{
\Theta(x_1^*x_2)\xi_2, \xi_1
}
\\
&\quad\quad
+
\ip{
\Theta(x_2^*x_1)\xi_1, \xi_2
}
+
\ip{
\Psi(x_2^*x_2)\xi_2, \xi_2
}.
\end{align*}
Therefore, 
\[
|s_\Theta(x_1 \otimes \xi_1, x_2 \otimes \xi_2)|^2 \leq s_\Phi(x_1 \otimes \xi_1, x_1 \otimes \xi_1) s_\Psi(x_2 \otimes \xi_2, x_2 \otimes \xi_2).
\]
Hence, by the maximality of the geometric mean, 
we have 
\[
s_\Theta \leq s_{\Phi \# \Psi},
\]
which implies
\[
\Theta \leq_{\mathrm{cp}} \Phi \# \Psi.
\]

Now we put $\Theta = \Phi \# \Psi$.
We next show that $s_\Pi$ is positive.
Let $x\in \cM$ and
\[
\boldsymbol{\xi} = 
\begin{bmatrix}
\eta_1 & \eta_2 \\
\zeta_1 & \zeta_2
\end{bmatrix}\in M_2(L^2(\cN)).
\]
By linearity and polarization, 
it suffices to check that
\[
s_\Pi(x \otimes \boldsymbol{\xi}, x \otimes \boldsymbol{\xi}) \geq 0.
\]
Note that
\begin{align*}
s_\Pi(x \otimes \boldsymbol{\xi}, x \otimes \boldsymbol{\xi})
&=
\sum_{i=1,2}\Big\{\ip{\Phi(x^*x)\eta_i, \eta_i} + \ip{\Theta(x^*x)\zeta_i, \eta_i} 
\\
&\qquad+ \ip{\Theta(x^*x)\eta_i, \zeta_i} 
+ \ip{\Psi(x^*x)\zeta_i, \zeta_i} \Big\}.
\end{align*}
Since $s_{\Phi \# \Psi} = s_\Phi \# s_\Psi$,
by the definition of the geometric mean, 
we have
\[
\ip{\Phi(x^*x)\eta_i, \eta_i} + \ip{\Theta(x^*x)\zeta_i, \eta_i} + \ip{\Theta(x^*x)\eta_i, \zeta_i} + \ip{\Psi(x^*x)\zeta_i, \zeta_i} \geq 0
\]
for $i=1,2$.
Therefore, $s_\Pi$ is positive, which means $\Pi$ is CP.
Thus the proposition follows.
\epf

\bthm\label{thm:main-geometric}
Let $\Phi, \Psi, \Phi_i, \Psi_i \in \mathrm{CP}(\cM, \cN)$ for $i=1,2$.
The following properties hold:
\benu
\item 
$\Phi_1 \leq_{\mathrm{cp}} \Phi_2$ and $\Psi_1 \leq_{\mathrm{cp}} \Psi_2$ 
imply 
$(\Phi_1 \# \Psi_1) \leq_{\mathrm{cp}} (\Phi_2 \# \Psi_2)$,
\item 
For any von Neumann algebra $\cP$ 
and $\Xi \in \mathrm{CP}(\cN, \cP)$,
\[
\Xi \circ (\Phi \# \Psi) \leq_{\mathrm{cp}} (\Xi \circ \Phi) \# (\Xi \circ \Psi),
\]
\item 
For any von Neumann algebra $\cL$ 
and $\Xi \in \mathrm{CP}(\cL, \cM)$,
\[
(\Phi \# \Psi) \circ \Xi \leq_{\mathrm{cp}} (\Phi \circ \Xi) \# (\Psi \circ \Xi),
\]
\item 
$\Phi_n, \Psi_n \in \mathrm{CP}(\cM, \cN)$, 
$\Phi_n \downarrow \Phi$ 
and $\Psi_n \downarrow \Psi$ 
imply $\Phi_n \# \Psi_n \downarrow \Phi \# \Psi$,
where $\Phi_n \downarrow \Phi$ means that 
$\Phi_1 \geq_{\mathrm{cp}} \Phi_2 \geq_{\mathrm{cp}} \cdots$ 
and $\Phi_n(x) \to \Phi(x)$ 
in the ultraweak operator topology for each $x \in \cM^+$.\eenu
\ethm

\bpf
(1) 
Let $\Theta = \Phi_1 \# \Psi_1$. 
By Proposition \ref{prop:geo}, 
we have
\[
\begin{bmatrix}
\Phi_1 & \Theta \\
\Theta & \Psi_1 \\
\end{bmatrix}
\geq_{\mathrm{cp}} 0.
\]
Since $\Phi_1 \leq_{\mathrm{cp}} \Phi_2$ and $\Psi_1 \leq_{\mathrm{cp}} \Psi_2$, 
we obtain
\[
\begin{bmatrix}
\Phi_2 & \Theta \\
\Theta & \Psi_2 \\
\end{bmatrix}
=
\begin{bmatrix}
\Phi_1 & \Theta \\
\Theta & \Psi_1 \\
\end{bmatrix}
+
\begin{bmatrix}
\Phi_2 - \Phi_1 & 0 \\
0 & \Psi_2 - \Psi_1 \\
\end{bmatrix}
\geq_{\mathrm{cp}} 0.
\]
By Proposition \ref{prop:geo}, 
this implies that 
$\Theta \leq_{\mathrm{cp}} \Phi_2 \# \Psi_2$, 
which is the desired conclusion.

(2) 
Let $\Xi \in \mathrm{CP}(\cN, \cP)$. 
By applying the CP map 
$\Xi^{(2)} = \mathrm{id}_{M_2(\mathbb{C})} \otimes \Xi$ to
\[
\begin{bmatrix}
\Phi & \Phi \# \Psi \\
\Phi \# \Psi & \Psi \\
\end{bmatrix}\geq_{\mathrm{cp}} 0,
\]
we obtain
\[
\begin{bmatrix}
\Xi \circ \Phi & \Xi \circ (\Phi \# \Psi) \\
\Xi \circ (\Phi \# \Psi) & \Xi \circ \Psi \\
\end{bmatrix}
=
\Xi^{(2)} \circ
\begin{bmatrix}
\Phi & \Phi \# \Psi \\
\Phi \# \Psi & \Psi \\
\end{bmatrix}
\geq_{\mathrm{cp}} 0,
\]
which implies 
\[
\Xi \circ (\Phi \# \Psi) \leq_{\mathrm{cp}} (\Xi \circ \Phi) \# (\Xi \circ \Psi).
\]

(3) 
Let $\Xi \in \mathrm{CP}(\cL, \cM)$. 
Since
\[
\begin{bmatrix}
\Phi \circ \Xi & (\Phi \# \Psi) \circ \Xi \\
(\Phi \# \Psi) \circ \Xi & \Psi \circ \Xi \\
\end{bmatrix}
=
\begin{bmatrix}
\Phi & \Phi \# \Psi \\
\Phi \# \Psi & \Psi \\
\end{bmatrix}
\circ \Xi
\geq_{\mathrm{cp}} 0,
\]
we obtain
\[
(\Phi \# \Psi) \circ \Xi \leq_{\mathrm{cp}} (\Phi \circ \Xi) \# (\Psi \circ \Xi).
\]

(4) 
Since $\Phi_n \geq_{\mathrm{cp}} \Phi_{n+1}$ 
and $\Psi_n \geq_{\mathrm{cp}} \Psi_{n+1}$, 
property (1) implies that
\[
\Phi_n \# \Psi_n \geq_{\mathrm{cp}} \Phi_{n+1} \# \Psi_{n+1}.
\]
Furthermore, since $\Phi_n \geq_{\mathrm{cp}} \Phi$ 
and $\Psi_n \geq_{\mathrm{cp}} \Psi$, 
we also have
\[
\Phi_n \# \Psi_n \geq_{\mathrm{cp}} \Phi \# \Psi \quad \text{for all } n.
\]
Thus, for each $x \in \cM^+$, 
the sequence
$\{(\Phi_n \# \Psi_n)(x)\}_{n=1}^\infty$ is monotonically decreasing  in $\cN^+$, 
and is bounded below by $(\Phi \# \Psi)(x)$.
Therefore, it converges in the ultraweak operator topology.
Since every element in $\cM$ is a linear combination of four positive elements, 
we can define a linear map $\Theta \colon \cM \to \cN$ 
by the point-ultraweak limit
\[
\Theta(x) \coloneqq \lim_{n \to \infty} (\Phi_n \# \Psi_n)(x)
\]
for $x \in \cM$.
Because $\mathrm{CP}(\cM,\cN)$ is closed under the point-ultraweak operator topology, 
$\Theta$ is CP.
By construction, 
we have 
\[
\Theta \geq_{\mathrm{cp}} \Phi \# \Psi.
\]
To show the reverse inequality, 
we use the block matrix characterization. 
For each $n$, 
we have
\[
\begin{bmatrix}
\Phi_n & \Phi_n \# \Psi_n \\
\Phi_n \# \Psi_n & \Psi_n \\
\end{bmatrix}
\geq_{\mathrm{cp}} 0.
\]
Taking the limit as $n \to \infty$ in the point-ultraweak operator topology, 
we obtain
\[
\begin{bmatrix}
\Phi & \Theta \\
\Theta & \Psi \\
\end{bmatrix}
\geq_{\mathrm{cp}} 0.
\]
By Proposition \ref{prop:geo}, 
this implies that $\Theta \leq_{\mathrm{cp}} \Phi \# \Psi$.
Therefore, we conclude that $\Theta = \Phi \# \Psi$, 
which means $\Phi_n \# \Psi_n \downarrow \Phi \# \Psi$.
\epf

\bcor\label{cor:geom-automorphism}
Let $\Phi, \Psi \in \mathrm{CP}(\cM, \cN)$. 
Suppose that $\alpha \in \mathrm{Aut}(\cN)$ and $\beta \in \mathrm{Aut}(\cM)$. 
Then 
\[
\alpha \circ (\Phi \# \Psi) \circ \beta = (\alpha \circ \Phi \circ \beta) \# (\alpha \circ \Psi \circ \beta).
\]
\ecor

\bpf
By Theorem \ref{thm:main-geometric}, we immediately have
\[
\alpha \circ (\Phi \# \Psi) \circ \beta \leq_{\mathrm{cp}} (\alpha \circ \Phi \circ \beta) \# (\alpha \circ \Psi \circ \beta).
\]
Applying Theorem \ref{thm:main-geometric} again, we obtain
\begin{align*}
&\alpha^{-1} \circ \left( (\alpha \circ \Phi \circ \beta) \# (\alpha \circ \Psi \circ \beta) \right) \circ \beta^{-1} 
\\
&\quad\leq_{\mathrm{cp}} (\alpha^{-1} \circ \alpha \circ \Phi \circ \beta \circ \beta^{-1}) \# (\alpha^{-1} \circ \alpha \circ \Psi \circ \beta \circ \beta^{-1}) \\
&\quad= \Phi \# \Psi.
\end{align*}
Applying $\alpha$ from the left and $\beta$ from the right to both sides of this inequality, we get
\[
(\alpha \circ \Phi \circ \beta) \# (\alpha \circ \Psi \circ \beta) \leq_{\mathrm{cp}} \alpha \circ (\Phi \# \Psi) \circ \beta.
\]
Combining these two inequalities yields the desired equality.
\epf

\bprop\label{prop:geo-concave}
Let $\Phi_i, \Psi_i\in \mathrm{CP}(\cM, \cN)$ for $i=1,2$.
Then
\[
(\Phi_1\#\Psi_1)+(\Phi_2\#\Psi_2)
\leq_{\mathrm{cp}}
(\Phi_1+\Phi_2)\#(\Psi_1+\Psi_2).
\]
\eprop

\bpf
Suppose that $\Theta_i\in \mathrm{CP}(\cM, \cN)$ satisfies
\[
\begin{bmatrix}
\Phi_i & \Theta_i \\
\Theta_i & \Psi_i \\
\end{bmatrix}\geq_{\mathrm{cp}} 0
\]
for $i=1,2$.
Then
\[
\begin{bmatrix}
\Phi_1+\Phi_2& \Theta_1+\Theta_2\\
\Theta_1+\Theta_2 & \Psi_1+\Psi_2\\
\end{bmatrix}\geq_{\mathrm{cp}} 0.
\]
By Proposition \ref{prop:geo}, 
\[
\Theta_1+\Theta_2\leq_{\mathrm{cp}}(\Phi_1+\Phi_2)\#(\Psi_1+\Psi_2).
\]
\epf

\subsection{Parallel sum}

We next consider the parallel sum for CP maps.
The following is essentially the same argument as for the geometric mean, but we provide a brief outline for the reader.

Let $\Phi, \Psi \in \mathrm{CP}(\cM, \cN)$.
By the Pusz--Woronowicz theory, we obtain the parallel sum $\gamma = s_\Phi : s_\Psi$ on $\cM \odot L^2(\cN)$.
We define a Hilbert space $\cH_\gamma$ 
via separation and completion
of $\cM\odot L^2(\cN)$ with respect to $\gamma$.
We denote by $h\colon \cM\odot L^2(\cN)\to\cH_\gamma$ the canonical map.

We first claim that the Hilbert space $\cH_\gamma$ has the left action of $\cM$ and the right action of $\cN$.
For a fixed $a\in\cM$, we define the sesquilinear form
$\gamma^a$ on $\cM\odot L^2(\cN)$ by
\[
\gamma^a(x\otimes \xi, y\otimes \eta)
\coloneqq
\gamma((ax)\otimes \xi, (ay)\otimes \eta)
\]
for $x, y\in\cM$ and $\xi, \eta\in L^2(\cN)$.
Suppose that $\zeta, \zeta_1, \zeta_2\in\cM\odot L^2(\cN)$
satisfy $\zeta=\zeta_1+\zeta_2$.
Then by Lemma \ref{lem:pw},
\[
\gamma(\zeta, \zeta)\leq s_\Phi(\zeta_1, \zeta_1)+s_\Psi(\zeta_2, \zeta_2).
\]
Since $a\zeta=a \zeta_1+a\zeta_2$, 
we have
\begin{align*}
\gamma^a(\zeta, \zeta)
&=
\gamma(a\zeta, a\zeta)
\\
&\leq
s_\Phi(a\zeta_1, a\zeta_1)+
s_\Psi(a\zeta_2, a\zeta_2)
\\
&\leq
\|a\|^2s_\Phi(\zeta_1, \zeta_1)
+\|a\|^2s_\Psi(\zeta_2, \zeta_2).
\end{align*}
By taking the infimum over $\zeta=\zeta_1+\zeta_2$,
we obtain
\[
\gamma^a\leq\|a\|^2(s_\Phi : s_\Psi).
\]
Therefore, we can define the left action of $\cM$ on $\cH_\gamma$ by
\[
a(x\otimes \xi)\coloneqq (ax)\otimes \xi
\]
for $x\in\cM$ and $\xi\in L^2(\cN)$.
Similarly, we can also define the right action of $\cN$ on $\cH_\gamma$ by
\[
(x\otimes \xi)b\coloneqq x\otimes (\xi b)
\]
for $x\in\cM$ and $\xi\in L^2(\cN)$.

Next we claim that an operator
\[
V\colon L^2(\cN)\to\cH_\gamma,
\quad
\xi\mapsto h(1\otimes \xi)
\]
is bounded.
Indeed, 
since $1\otimes \xi=1\otimes \xi+0=0+1\otimes \xi$,
we have
\begin{align*}
\gamma(1\otimes \xi, 1\otimes \xi)
&\leq
s_\Phi(1\otimes \xi, 1\otimes \xi)+
s_\Psi(0, 0)
\\
&\leq
\|\Phi\|\|\xi\|^2.
\end{align*}
Similarly,
\[
\gamma(1\otimes \xi, 1\otimes \xi)
\leq
\|\Psi\|\|\xi\|^2.
\]
Hence,
\[
\gamma(1\otimes \xi, 1\otimes \xi)
\leq\min\{\|\Phi\|, \|\Psi\|\}\|\xi\|^2,
\]
which implies that $V$ is bounded.
Then we define the CP map $\Theta$ on $\cM$ by
\[
\Theta(x)=V^*xV
\]
for $x\in\cM$.

We finally show that $\Theta(x)\in\cN$ for $x\in\cM$.
By (\ref{eq:right-action}), for $b\in\cN$, 
we have
\[
\ip{
\Theta(x)Jb^*J\xi, \eta
}
=
\ip{
Jb^*J\Theta(x)\xi, \eta
}.
\]
Therefore $\Theta(x)\in(\cN')'=\cN$.
By (\ref{eq:left-action}),
we also have
\[
\ip{
\Theta(y^*x)\xi, \eta
}
=
(s_\Phi : s_\Psi)(x\otimes \xi, y\otimes \eta).
\]

\bdf
For $\Phi, \Psi\in\mathrm{CP}(\cM, \cN)$,
the {\em parallel sum} $\Phi : \Psi$
is defined by
\[
\Phi : \Psi \coloneqq \Theta.
\]
\edf

Therefore, the norm estimate and other basic properties are summarized as follows.

\bprop\label{prop:parallel-basic}
Let $\Phi, \Psi \in \mathrm{CP}(\cM, \cN)$.
Then the following properties hold:
\benu
\item $\Phi : \Psi = \Psi : \Phi$,
\item $\Phi : \Phi = \frac{1}{2}\Phi$,
\item $\Phi : \Psi\leq_{\mathrm{cp}}\Phi$ and $\Phi : \Psi\leq_{\mathrm{cp}}\Psi$.
In particular,
$\|\Phi : \Psi\| \leq \min\{\|\Phi\|, \|\Psi\|\}$.
\eenu
\eprop

Before stating the fundamental properties of the parallel sum, we provide its block matrix characterization. 
As with the geometric mean, the corresponding characterization for CP maps follows directly from Proposition \ref{prop:harmonic} (1).

\bprop\label{prop:parallel-matrix}
Let $\Phi, \Psi \in \mathrm{CP}(\cM, \cN)$. 
Then
\[
\Phi : \Psi = \max
\left\{
\Theta \in \mathrm{CP}(\cM, \cN) \;\middle|\;
\begin{bmatrix}
\Phi & \Phi \\
\Phi & \Phi + \Psi \\
\end{bmatrix}
\geq_{\mathrm{cp}}
\begin{bmatrix}
\Theta & 0 \\
0 & 0 \\
\end{bmatrix}
\right\}.
\]
\eprop

By employing the same arguments as in the case of the geometric mean,
we can establish analogous fundamental properties for the parallel sum.

\bthm\label{thm:main-parallel}
Let $\Phi, \Psi, \Phi_i, \Psi_i \in \mathrm{CP}(\cM, \cN)$ for $i=1,2$.
The following properties hold:
\benu
\item $\Phi_1 \leq_{\mathrm{cp}} \Phi_2$ and $\Psi_1 \leq_{\mathrm{cp}} \Psi_2$ imply $(\Phi_1 : \Psi_1) \leq_{\mathrm{cp}} (\Phi_2 : \Psi_2)$,
\item For any von Neumann algebra $\cP$ and $\Xi \in \mathrm{CP}(\cN, \cP)$,
\[
\Xi \circ (\Phi : \Psi) \leq_{\mathrm{cp}} (\Xi \circ \Phi) : (\Xi \circ \Psi),
\]
\item For any von Neumann algebra $\cL$ and $\Xi \in \mathrm{CP}(\cL, \cM)$,
\[
(\Phi : \Psi) \circ \Xi \leq_{\mathrm{cp}} (\Phi \circ \Xi) : (\Psi \circ \Xi),
\]
\item $\Phi_n \downarrow \Phi$ and $\Psi_n \downarrow \Psi$ imply $\Phi_n : \Psi_n \downarrow \Phi : \Psi$.
\eenu
\ethm

\bcor\label{cor:para-automorphism}
Let $\Phi_1, \Phi_2 \in \mathrm{CP}(\cM, \cN)$. 
Suppose that $\alpha \in \mathrm{Aut}(\cN)$ and $\beta \in \mathrm{Aut}(\cM)$. 
Then 
\[
\alpha \circ (\Phi_1 : \Phi_2) \circ \beta = (\alpha \circ \Phi_1 \circ \beta) : (\alpha \circ \Phi_2 \circ \beta).
\]
\ecor

\bprop\label{prop:para-concave}
Let $\Phi_i, \Psi_i\in \mathrm{CP}(\cM, \cN)$ for $i=1,2$.
Then
\[
(\Phi_1 : \Psi_1)+(\Phi_2 : \Psi_2)
\leq_{\mathrm{cp}}
(\Phi_1+\Phi_2) : (\Psi_1+\Psi_2).
\]
\eprop

Finally, to bridge to the next section, we introduce the harmonic mean for CP maps. 
In the Kubo--Ando theory, the harmonic mean is simply a scalar multiple of the parallel sum. We adopt the same definition here.

\bdf\label{def:harmonic}
For $\Phi, \Psi \in \mathrm{CP}(\cM, \cN)$, the \emph{harmonic mean} $\Phi\, !\, \Psi$ is defined by
\[
\Phi\, !\, \Psi \coloneqq 2 (\Phi : \Psi).
\]
\edf

The block matrix characterization
of the harmonic mean also follows from Proposition \ref{prop:harmonic} (2).

\bprop\label{prop:harmonic-cp}
Let $\Phi, \Psi \in \mathrm{CP}(\cM, \cN)$. Then
\[
\Phi\, !\, \Psi = \max
\left\{
\Theta \in \mathrm{CP}(\cM, \cN) \;\middle|\;
\begin{bmatrix}
2\Phi & 0 \\
0 & 2\Psi \\
\end{bmatrix}
\geq_{\mathrm{cp}}
\begin{bmatrix}
\Theta & \Theta \\
\Theta & \Theta \\
\end{bmatrix}
\right\}.
\]
\eprop


\subsection{AM-GM-HM inequality}

In this section, we establish the Arithmetic-Geometric-Harmonic Mean (AM-GM-HM) inequality for CP maps.

\bdf
Let $\Phi, \Psi \in \mathrm{CP}(\cM, \cN)$.
The {\em arithmetic mean} $\Phi\nabla\Psi$ 
is naturally defined by 
\[
\Phi\nabla\Psi\coloneqq\frac{\Phi + \Psi}{2}.
\]
\edf

This approach highlights the intrinsic order-theoretic nature of the
AM--GM--HM inequality and provides a perspective compatible with the
connection framework developed later.

\bthm[AM-GM-HM inequality]\label{thm:am-gm-hm}
Let $\Phi, \Psi \in \mathrm{CP}(\cM, \cN)$. Then
\[
\Phi\, !\, \Psi \leq_{\mathrm{cp}} \Phi \# \Psi \leq_{\mathrm{cp}} \Phi\nabla\Psi.
\]
\ethm

\bpf
By {\cite[Theorem 1.1]{pw1}}, 
the associated forms $s_\Phi$ and $s_\Psi$ 
can be represented 
by compatible representations,
\[
s_\Phi(\zeta_1, \zeta_2) = \ip{Ah(\zeta_1), h(\zeta_2)} \quad \text{and} \quad s_\Psi(\zeta_1, \zeta_2) = \ip{B h(\zeta_1), h(\zeta_2)}.
\]
By the definition of the operator means for positive sesquilinear forms, 
the forms $s_{\Phi\, !\, \Psi}$, $s_{\Phi \# \Psi}$, and $s_{\Phi \nabla \Psi}$ correspond exactly to the operators $A\, !\, B$, $A \# B$, and $A\nabla B$, respectively.
The Kubo--Ando theory guarantees that
\[
A\, !\, B \leq A \# B \leq A\nabla B.
\]
Taking the inner product with $h(\zeta)$
immediately yields that
\[
s_{\Phi\, !\, \Psi}\leq s_{\Phi \# \Psi} \leq s_{\Phi \nabla \Psi}.
\]
Since the correspondence $\Theta \mapsto s_\Theta$ preserves the order, we finally conclude that
\[
\Phi\, !\, \Psi \leq_{\mathrm{cp}} \Phi \# \Psi \leq_{\mathrm{cp}} \Phi \nabla \Psi.
\]
\epf

\bcor\label{cor:am-gm-equality}
Let $\Phi, \Psi \in \mathrm{CP}(\cM, \cN)$. 
Then the following conditions are equivalent:
\benu
\item $\Phi \nabla \Psi = \Phi \# \Psi$,
\item $\Phi \# \Psi = \Phi\, !\, \Psi$,
\item $\Phi = \Psi$.
\eenu
\ecor

\bpf
It is trivial that (3) implies both (1) and (2).

To prove that (1) implies (3), suppose that $\Phi \nabla \Psi = \Phi \# \Psi$. 
The corresponding bounded positive operators $A, B$ constructed in the proof of Theorem \ref{thm:am-gm-hm} must satisfy
\[
A\nabla B = A \# B.
\]
By the Kubo--Ando theory, 
the equality condition for the AM-GM inequality holds if and only if $A = B$. 
Thus, we have $A = B$, which implies that 
$s_\Phi= s_\Psi$. 
Hence, we conclude that
$\Phi = \Psi$. 

The implication from (2) to (3) follows by exactly the same argument using the equality condition for the GM-HM inequality.
\epf

\bcor\label{cor:unital-bound}
Let $\Phi, \Psi \in \mathrm{CP}(\cM, \cN)$ be unital.
Then 
\[
(\Phi\, !\, \Psi)(1) \leq (\Phi \# \Psi)(1) \leq 1.
\]
\ecor

\bpf
Evaluating the AM-GM-HM inequality at the identity $1 \in \cM$, we obtain
\[(\Phi\, !\, \Psi)(1) \leq (\Phi \# \Psi)(1) \leq (\Phi \nabla \Psi)(1) = \frac{\Phi(1) + \Psi(1)}{2} = \frac{1 + 1}{2} = 1.\]
\epf

\section{Examples in matrix algebras}

In this section, we consider the finite-dimensional case to illustrate our main results. 
Let $\cM = M_m(\mathbb{C})$ 
and $\cN = M_n(\mathbb{C})$. 
We denote by $\Tr$ the canonical unnormalized trace.
Recall that every CP map $\Phi \colon M_m(\mathbb{C}) \to M_n(\mathbb{C})$ 
is uniquely represented by its Choi matrix $C_\Phi \in M_m(\mathbb{C}) \otimes M_n(\mathbb{C})$ defined by
\[
C_\Phi \coloneqq \sum_{i,j=1}^m e_{ij} \otimes \Phi(e_{ij}),
\]
where $\{e_{ij}\}_{i,j=1}^m$ are the standard matrix units of $M_m(\mathbb{C})$.
This correspondence is called the Choi--Jamio\l{}kowski isomorphism.
The map $\Phi \mapsto C_\Phi$ is an affine isomorphism that preserves the order, i.e., $\Phi \leq_{\mathrm{cp}} \Psi$ if and only if $C_\Phi \leq C_\Psi$ as positive semidefinite matrices. 
Using this correspondence, the geometric mean of CP maps naturally reduces to the geometric mean of their Choi matrices.

\bprop\label{prop:choi-matrix}
Let $\Phi, \Psi \colon M_m(\mathbb{C}) \to M_n(\mathbb{C})$ be CP maps, 
and let $C_\Phi, C_\Psi \in M_m(\mathbb{C}) \otimes M_n(\mathbb{C})$ be their Choi matrices. Then, the Choi matrix $C_{\Phi \# \Psi}$ of the geometric mean $\Phi \# \Psi$ is exactly the geometric mean of $C_\Phi$ and $C_\Psi$, that is,
\[
C_{\Phi \# \Psi} = C_\Phi \# C_\Psi.
\]
\eprop

\bpf
Let $\Theta \colon M_m(\mathbb{C}) \to M_n(\mathbb{C})$ be a CP map. 
By Proposition \ref{prop:geo}, $\Theta \leq_{\mathrm{cp}} \Phi \# \Psi$ if and only if the block map
\[
\Pi = \begin{bmatrix}
\Phi & \Theta \\
\Theta & \Psi \\
\end{bmatrix}
\colon M_m(\mathbb{C}) \to M_2(M_n(\mathbb{C}))
\]
is CP. 
Then a direct computation of the corresponding Choi matrix shows that
\[
C_\Pi = \begin{bmatrix}
C_\Phi & C_\Theta \\
C_\Theta & C_\Psi \\
\end{bmatrix}.
\]
Since $\Pi$ is CP if and only if $C_\Pi \geq 0$, we have that $\Theta \leq_{\mathrm{cp}} \Phi \# \Psi$ is equivalent to
\[
\begin{bmatrix}
C_\Phi & C_\Theta \\
C_\Theta & C_\Psi \\
\end{bmatrix}
\geq 0.
\]
Recall that the geometric mean $C_\Phi \# C_\Psi$ is defined as the maximum positive matrix $X$ satisfying 
\[
\begin{bmatrix} 
C_\Phi & X \\ 
X & C_\Psi 
\end{bmatrix} \geq 0.
\]
Because the Choi--Jamio\l{}kowski isomorphism $\Theta \mapsto C_\Theta$ is order-preserving, the maximal CP map $\Phi \# \Psi$ corresponds exactly to the maximal matrix $C_\Phi \# C_\Psi$. 
Therefore, we conclude that $C_{\Phi \# \Psi} = C_\Phi \# C_\Psi$.
\epf

\brem
By using the block matrix characterization for the harmonic mean (Proposition \ref{prop:harmonic-cp}) and applying exactly the same argument as in Proposition \ref{prop:choi-matrix}, we can easily deduce that the Choi matrix of the harmonic mean is the harmonic mean of the respective Choi matrices:
\[
C_{\Phi\, !\, \Psi} = C_\Phi\, !\, C_\Psi.
\]
By the Choi--Jamio\l{}kowski isomorphism, 
the AM-GM-HM inequality for CP maps (Theorem \ref{thm:am-gm-hm}) in the finite-dimensional case perfectly translates to the classical AM-GM-HM inequality for positive matrices:
\[
C_\Phi\, !\, C_\Psi \leq C_\Phi \# C_\Psi \leq C_\Phi \nabla C_\Psi.
\]
This correspondence ensures that our generalizations of the operator means to CP maps are fully consistent with the classical matrix theory via the Choi--Jamio\l{}kowski isomorphism.
\erem

\brem\label{rem:projections}
Let $\Phi,\Psi \colon M_m(\mathbb C)\to M_n(\mathbb C)$ be CP maps with
\[
C_\Phi = rP,
\qquad
C_\Psi = sQ
\]
for some projections $P, Q$ and scalars $r,s\geq 0$.
Then 
\[
C_\Phi\#C_\Psi=\sqrt{rs}(P\wedge Q)
\quad\text{and}\quad
C_\Phi\, !\, C_\Psi=\frac{2rs}{r+s}(P\wedge Q).
\]

Indeed, 
for projections $P,Q$, it is well-known that
$P\# Q=P\, !\, Q=P\wedge Q$.
Therefore, we obtain
\[
C_\Phi\# C_\Psi=(rP)\#(sQ)
=\sqrt{rs}(P\# Q)
=\sqrt{rs}(P\wedge Q),
\]
and
\[
C_\Phi\, !\, C_\Psi
=(rP)\, !\, (sQ)
=\frac{2rs}{r+s}(P\, !\, Q)
=\frac{2rs}{r+s}(P\wedge Q).
\]
\erem

\bex[Non-unitality]\label{ex:non-unital}
We note that even if $\Phi$ and $\Psi$ are unital CP maps, their geometric mean $\Phi \# \Psi$ need not be unital.

Let $\mathcal M = \mathcal N = M_2(\mathbb C)$.
Consider the identity map $\Phi(x)=x$ and the unitary conjugation
\[
\Psi(x)=UxU^*,\qquad
U=
\begin{bmatrix}
1 & 0 \\
0 & -1  \\
\end{bmatrix}.
\]
Clearly, $\Phi(1)=1$ and $\Psi(1)=1$.
Their Choi matrices $C_\Phi$, $C_\Psi$ are
\[
C_\Phi=\sum_{i,j=1}^2 e_{ij}\otimes e_{ij},
\quad
C_\Psi=\sum_{i,j=1}^2 e_{ij}\otimes (Ue_{ij}U^*).
\]
Let $\{e_i\}_{i=1,2}$ be the canonical orthonormal basis of $\mathbb C^2$. 
Note that $e_{ij}=e_ie_j^*$.
Define
\[
v_\Phi=e_1\otimes e_1+e_2\otimes e_2,
\quad
v_\Psi=e_1\otimes e_1-e_2\otimes e_2
\in\mathbb C^2\otimes\mathbb C^2 .
\]
Then
\[
C_\Phi=v_\Phi v_\Phi^*,
\qquad
C_\Psi=v_\Psi v_\Psi^*,
\]
and so both Choi matrices are positive rank-one operators.
Let $P$ and  $Q$ be the projections onto $\C v_\Phi$ and $\C v_\Psi$, respectively.
Then
\[
C_\Phi=\|v_\Phi\|^2 P,
\qquad
C_\Psi=\|v_\Psi\|^2 Q.
\]
Since $\langle v_\Phi, v_\Psi \rangle = 0$, we have $v_\Phi\perp v_\Psi$.
Therefore,
$\mathrm{ran}(P)\cap\mathrm{ran}(Q)=\{0\}$.
By Remark \ref{rem:projections}, we obtain
$\Phi\#\Psi=0$.
\eex

\bex[Density matrices]\label{ex:states}
Let $\cM=M_n(\mathbb C)$ and $\cN=\mathbb C$.
A positive functional on $M_n(\mathbb C)$ is a CP map
$\varphi \colon M_n(\mathbb C)\to\mathbb C$.
Every positive functional $\varphi$ is represented by a density matrix
$\rho\in M_n(\mathbb C)^+$ such that for $x\in M_n(\mathbb C)$,
\[
\varphi(x)=\mathrm{Tr}(\rho x).
\]
Similarly, let $\psi$ be another positive functional represented by a density matrix
$\sigma\in M_n(\mathbb C)^+$.

We compute the Choi matrix of $\varphi$.
Since $M_n(\mathbb C)\otimes\mathbb C$ is naturally identified with $M_n(\mathbb C)$,
we obtain
\[
C_\varphi
=
\sum_{i,j=1}^n e_{ij}\otimes\varphi(e_{ij})
=
\sum_{i,j=1}^n \mathrm{Tr}(\rho e_{ij}) e_{ij}
=
\sum_{i,j=1}^n \rho_{ji} e_{ij}
=
\rho^T .
\]
Thus, the Choi matrix of a positive functional is the transpose of its density matrix.
Similarly, we have $C_\psi=\sigma^T$.

Now we consider the geometric mean $\varphi \#\psi$.
Since the transpose map $A\mapsto A^T$ is an order-preserving isomorphism on
positive semidefinite matrices, it commutes with operator means.
Hence,
\[
C_{\varphi \#\psi}
=
C_\varphi \# C_\psi
=
\rho^T \#\sigma^T
=
(\rho \#\sigma)^T .
\]
By the Choi--Jamio\l{}kowski isomorphism, 
we conclude that
\[
(\varphi \#\psi)(x)
=\mathrm{Tr}((\rho \#\sigma)x).
\]
Therefore, the geometric mean of positive functionals corresponds exactly to
the geometric mean of their density matrices.
\eex

\bex[Identity channel and completely depolarizing channel]\label{ex:quantum-channels}

In quantum information theory, two fundamental extremal channels are the identity channel, representing perfect transmission of quantum information, and the completely depolarizing channel, representing maximal noise.
We compute the geometric and harmonic means of these two channels.

Let $\mathcal{M}=\mathcal{N}=M_d(\mathbb{C})$ with $d\ge2$.
Consider the identity channel
\[
\mathrm{id}(x)=x
\]
and the completely depolarizing channel
\[
\Delta(x)=\frac{\mathrm{Tr}(x)}{d}\, 1 .
\]
Both maps are unital CP.

Let $C_{\mathrm{id}},C_{\Delta}\in M_d(\mathbb{C})\otimes M_d(\mathbb{C})$ denote their Choi matrices.
They are given by
\[
C_{\mathrm{id}}=\sum_{i,j=1}^d e_{ij}\otimes e_{ij},
\qquad
C_{\Delta}
=\sum_{i,j=1}^d e_{ij}\otimes\Delta(e_{ij})
=\sum_{i=1}^d e_{ii}\otimes\frac1d1
=\frac1d(1\otimes1).
\]
Let
\[
v=\sum_{i=1}^d e_i\otimes e_i\in\mathbb{C}^d\otimes\mathbb{C}^d .
\]
Then
\[
C_{\mathrm{id}}=vv^*
\]
is a rank-one positive operator.
In contrast,
\[
C_\Delta=\frac{1}{d} I
\]
is a scalar multiple of the identity on $\mathbb{C}^d\otimes\mathbb{C}^d$.
Let $P$ be the projection onto $\C v$,
and $Q=I$.
Then
\[
C_{\mathrm{id}}= dP,\qquad C_\Delta=\frac{1}{d}Q.
\]
Since $\mathrm{ran}(P)\cap\mathrm{ran}(Q)=\mathrm{ran}(P)$, by Remark \ref{rem:projections},
we have
\[
C_{\mathrm{id}}\# C_\Delta
=\sqrt{d\cdot\frac{1}{d}}(P\wedge Q)=P=\frac{1}{d}C_{\mathrm{id}},
\]
and
\[
C_{\mathrm{id}}\, !\, C_\Delta
=\frac{2dd^{-1}}{d+d^{-1}}(P\wedge Q)=\frac{2d}{d^2+1}P=\frac{2}{d^2+1}C_{\mathrm{id}}.
\]
Thus the geometric mean corresponds to the projection onto the maximally entangled state. 

By the Choi--Jamio\l{}kowski isomorphism, the resulting quantum channels satisfy
\[
\mathrm{id} \# \Delta=\frac1d\, \mathrm{id},
\qquad
\mathrm{id}\, !\, \Delta=\frac{2}{d^2+1}\, \mathrm{id}.
\]
\eex

\bex[Schur multipliers]\label{ex:schur-multiplier}
Let $\cM=\cN=M_n(\mathbb C)$.
For $A\in M_n(\mathbb C)$, the Schur multiplier
$S_A\colon M_n(\mathbb C)\to M_n(\mathbb C)$ is defined by
\[
S_A(x)=A\circ x,
\]
where $\circ$ denotes the entrywise (Hadamard) product.
It is well-known that $S_A$ is CP if and only if $A$ is positive semidefinite.

We compute the Choi matrix of $S_A$.
Using the canonical matrix units $\{e_{ij}\}_{i,j=1}^n$, 
we obtain
\begin{align*}
C_{S_A}
&=
\sum_{i,j=1}^n e_{ij}\otimes S_A(e_{ij})
\\
&=
\sum_{i,j=1}^n e_{ij}\otimes (A_{ij}e_{ij})
\\
&=\sum_{i,j=1}^n A_{ij}(e_i e_j^*\otimes e_i e_j^*).
\end{align*}
Observe that
\[
e_i e_j^*\otimes e_i e_j^*
=
(e_i\otimes e_i)(e_j\otimes e_j)^* .
\]
Let $V \colon \mathbb C^n\to\mathbb C^n\otimes\mathbb C^n$ be the isometry defined by
\[
V e_i = e_i\otimes e_i
\]
for $i=1,\dots,n$.
Then the above expression can be written as
\[
C_{S_A}=VAV^* .
\]

Now let $A, B\ge 0$.
Since $V$ is an isometry, the transformer inequality for operator means yields
\[
(VAV^*) \# (VBV^*) = V(A \# B)V^* .
\]
Therefore,
\[
C_{S_A}\# C_{S_B}
=
V(A\#B)V^* .
\]
Recognizing that $V(A\#B)V^*$ is precisely the Choi matrix of the Schur multiplier
$S_{A\#B}$, we conclude that
\[
S_A\# S_B = S_{A\#B}.
\]
\eex

\bex[Adjoint maps]
Let $\cM = \cN = M_2(\mathbb{C})$.
For $A \in M_2(\mathbb{C})$, we define a CP map
\[
\Psi_A(x) = A x A^*
\]
for $x \in M_2(\mathbb{C})$.
Let
\[
A =
\begin{bmatrix}
2 & 0 \\
0 & 1 \\
\end{bmatrix},
\quad
B =
\begin{bmatrix}
1 & 0 \\
0 & 2 \\
\end{bmatrix}.
\]
Since $A$ and $B$ commute, their geometric mean is
\[
C = A \# B = A^{1/2} B^{1/2}
=
\begin{bmatrix}
\sqrt{2} & 0 \\
0 & \sqrt{2} \\
\end{bmatrix}
= \sqrt{2}\,I.
\]
Hence
\[
\Psi_C(x) = (\sqrt{2}I)x(\sqrt{2}I)^* = 2x.
\]

We now compute the geometric mean $\Psi_A\# \Psi_B$.
Their Choi matrices are rank-one operators
\[
C_{\Psi_A}=v_Av_A^*,
\qquad
C_{\Psi_B}=v_Bv_B^*,
\]
where
\[
v_A =
\begin{bmatrix}
2 \\ 0 \\ 0 \\ 1 \\
\end{bmatrix},
\qquad
v_B =
\begin{bmatrix}
1 \\ 0 \\ 0 \\ 2 \\
\end{bmatrix}.
\]
Let $P$ and $Q$ be the projections onto $\C v_A$ and $\C v_B$, respectively.
Then we have
\[
C_{\Psi_A}=\|v_A\|^2 P,
\qquad
C_{\Psi_B}=\|v_B\|^2 Q.
\]
Since $v_A, v_B$ are linearly independent,
$\mathrm{ran}(P)\cap\mathrm{ran}(Q)=\{0\}$.
By Remark \ref{rem:projections},
we obtain
\[
C_{\Psi_A} \# C_{\Psi_B}=(\|v_A\|^2 P)\#(\|v_B\|^2 Q)
=\|v_A\|\|v_B\|(P\wedge Q)=0.
\]
Consequently,
\[
\Psi_A \# \Psi_B = 0.
\]
Therefore, in general,
\[
\Psi_A \# \Psi_B \neq \Psi_{A\# B}.
\]
Thus, for adjoint CP maps, the geometric mean does not behave
functorially with respect to the operator geometric mean of the
implementing operators.
\eex

\brem\label{rem:kraus-support}
For general CP maps $\Phi$ and $\Psi$ with Kraus representations 
\[
\Phi(x) = \sum_{i} A_i x A_i^*
\quad\text{and}\quad 
\Psi(x) = \sum_{j} B_j x B_j^*,
\]
it is nontrivial to compute an explicit Kraus representation for their geometric mean $\Phi \# \Psi$. 
Their Choi matrices are given by 
\[
C_\Phi = \sum_{i} v_{A_i} v_{A_i}^*
\quad\text{and}\quad
C_\Psi = \sum_{j} v_{B_j} v_{B_j}^*.
\]
Since 
\[
C_\Phi \# C_\Psi = C_\Phi^{1/2}(C_\Phi^{-1/2}C_\Psi C_\Phi^{-1/2})^{1/2}C_\Phi^{1/2},
\]
if $C_\Phi$ and $C_\Psi$ are invertible, there is no simple formula to express the Kraus operators of $\Phi \# \Psi$ in terms of $A_i$ and $B_j$.

However, the spatial structure of the geometric mean is 
described by operator mean theory. 
As shown by \cite{fujii1992}, if $\mathrm{ran}(A)$ and $\mathrm{ran}(B)$ 
are closed, then
\[
\mathrm{ran}(A \# B)
=
\mathrm{ran}(A) \cap \mathrm{ran}(B).
\]
Hence we obtain
\[
\mathrm{ran}(C_{\Phi \# \Psi})
=
\mathrm{span}\{v_{A_i}\}_i
\cap
\mathrm{span}\{v_{B_j}\}_j.
\]
\erem

\section{Geometric mean in von Neumann algebras}

\subsection{Geometric mean of normal positive functionals}

Let $(\cM, \cH, J, \cP)$ be a standard form of a von Neumann algebra $\cM$. 
Let $\varphi$ and $\psi$ be normal positive functionals on $\cM$. 
To define their geometric mean, we consider the dominating functional
\[
\omega \coloneqq \varphi + \psi .
\]
We denote its unique implementing vector by
$\Omega \in \cP$.
Since $\varphi \leq \omega$ and $\psi \leq \omega$, 
there exist unique positive operators $h', k' \in \cM'$ such that
\[
\varphi(x) = \langle x h' \Omega, \Omega \rangle,
\qquad
\psi(x) = \langle x k' \Omega, \Omega \rangle
\]
for $x \in \cM$.
Note that $h'+k'=1_{\cH}$,
in particular, $h'$ and $k'$ commute.

Thus the Pusz--Woronowicz geometric mean $\varphi \# \psi$ is given by
\[
(\varphi \# \psi)(x)
=
\langle x (h' \# k') \Omega, \Omega \rangle
=
\langle x \sqrt{h'k'}\, \Omega, \Omega \rangle
\]
for $x \in \cM$.

\brem\label{rem:kosaki-vs-pw}
It is instructive to compare the Pusz--Woronowicz mean with 
another interpolation of states introduced in \cite{am, k1}, 
which is based on Connes' relative modular operator 
$\Delta_{\varphi,\psi}$.

We show that these two constructions produce different states in the 
noncommutative setting by computing them explicitly in the matrix case 
$\cM=M_n(\mathbb C)$.

Let $\varphi(x)=\mathrm{Tr}(\rho x)$ and $\psi(x)=\mathrm{Tr}(\sigma x)$ 
for invertible positive semidefinite matrices $\rho$ and $\sigma$, respectively.
We use the standard Hilbert space 
$\cH=M_n(\mathbb C)$ equipped with the Hilbert--Schmidt inner product
\[
\langle x,y\rangle=\mathrm{Tr}(y^*x).
\]
The representative vectors in $\cP$ are
\[
\Omega_\varphi=\rho^{1/2}, 
\qquad 
\Omega_\psi=\sigma^{1/2}.
\]
The relative modular operator is given by
\[
\Delta_{\varphi,\psi}(x)=\rho x\sigma^{-1},
\]
and
\[
\Delta_{\varphi,\psi}^{1/4}(x)
=
\rho^{1/4}x\sigma^{-1/4}
\quad\text{for}\ x\in\cH.
\]
In particular,
\[
\Delta_{\varphi,\psi}^{1/4}\Omega_\psi
=
\rho^{1/4}\sigma^{1/4}.
\]

The resulting state is
\[
\omega_{1/2}(x)
=
\mathrm{Tr}\!\left(
(\rho^{1/4}\sigma^{1/4})^*
x
(\rho^{1/4}\sigma^{1/4})
\right)
=
\mathrm{Tr}(\rho^{1/4}\sigma^{1/2}\rho^{1/4}x).
\]
In particular,
\[
\omega_{1/2}(1)
=
\mathrm{Tr}(\rho^{1/2}\sigma^{1/2}).
\]
On the other hand,
\[
\mathrm{Tr}(\rho\#\sigma)
=
\mathrm{Tr}\!\left(
\rho^{1/2}
(\rho^{-1/2}\sigma\rho^{-1/2})^{1/2}
\rho^{1/2}
\right).
\]

These quantities are naturally related to the 
Uhlmann fidelity in quantum information theory,
defined by
\[
F(\rho,\sigma)
=
\mathrm{Tr}\!\sqrt{\rho^{1/2}\sigma\rho^{1/2}}.
\]

It is known that
\[
\mathrm{Tr}(\rho\#\sigma)
\le
\mathrm{Tr}(\rho^{1/2}\sigma^{1/2})
\le
\mathrm{Tr}\sqrt{\rho^{1/2}\sigma\rho^{1/2}},
\]
and equality holds if and only if $\rho$ and $\sigma$ commute 
(see \cite{ara,h-p}). 
\erem

\subsection{Geometric mean of conditional expectations}

Conditional expectations are fundamental objects in the theory of 
von Neumann algebras, particularly in subfactor theory and 
quantum probability. 
It is therefore natural to investigate how operator means interact 
with conditional expectations.

Let $\cN \subseteq \cM$ be an inclusion of von Neumann algebras.
Recall that a projection from $\cM$ onto $\cN$ is a linear map $E \colon \cM \to \cN$
such that $E(a)=a$ for every $a\in\cN$.
A conditional expectation from $\cM$ onto $\cN$ is a contractive CP projection
$E \colon \cM \to \cN$ such that
$E(a x b)=aE(x)b$
for $a$, $b\in \cN$, $x\in \cM$.
A normal conditional expectation preserving a fixed faithful normal state is unique by 
Takesaki's theorem \cite{tak}.

Let 
\[
E_i \colon \cM \to \cN_i 
\]
be conditional expectations for $i=1,2$.
Composing with the natural inclusions $j_i \colon \cN_i \hookrightarrow \cM$, 
we regard
\[
\tilde E_i \coloneqq j_i \circ E_i
\]
as CP maps on $\cM$.

\bprop\label{prop:ce-bimodule}
Let $E_i \colon \cM\to\cN_i$ be conditional expectations for $i=1,2$ and set
\[
\Theta \coloneqq \tilde E_1 \# \tilde E_2 \in \mathrm{CP}(\cM,\cM).
\]
Then the following hold:

\begin{enumerate}
\item $\Theta$ is an $(\cN_1\cap \cN_2)$-bimodule map, i.e.,
\[
\Theta(a x b)=a\Theta(x)b
\]
for $a,b\in \cN_1\cap\cN_2,\ x\in \cM$,
\item $\Theta$ is subunital and
\[
\Theta(1)\in (\cN_1\cap\cN_2)' \cap \cM .
\]
\end{enumerate}
\eprop

\bpf

(1)
We first show that $\Theta$ intertwines inner automorphisms 
implemented by elements in $\cN_1 \cap \cN_2$, 
and then extend the result by linearity.
Let $a\in \cN_1\cap\cN_2$ be invertible. 
We define
$\Psi_a(x)=axa^*$.

Since
\[
\begin{bmatrix}
\tilde E_1 & \Theta \\
\Theta & \tilde E_2
\end{bmatrix}
\ge_{\mathrm{cp}}0,
\]
we have
\[
\begin{bmatrix}
\Psi_a\circ \tilde E_1 & \Psi_a\circ\Theta\\
\Psi_a\circ\Theta & \Psi_a\circ \tilde E_2
\end{bmatrix}
\ge_{\mathrm{cp}}0.
\]
Since $a \in \cN_i$ and $E_i$ is a conditional expectation onto $\cN_i$,
we have $E_i(axa^*) = a E_i(x) a^*$, i.e.,
\[
\tilde E_i \circ \Psi_a = \Psi_a \circ \tilde E_i.
\]
Then we obtain
\[
\begin{bmatrix}
\tilde E_1\circ\Psi_a & \Psi_a\circ\Theta\\
\Psi_a\circ\Theta & \tilde E_2\circ\Psi_a
\end{bmatrix}
\ge_{\mathrm{cp}}0.
\]
Therefore,
\[
\begin{bmatrix}
\tilde E_1 & \Psi_a\circ\Theta\circ\Psi_{a^{-1}}\\
\Psi_a\circ\Theta\circ\Psi_{a^{-1}} & \tilde E_2
\end{bmatrix}
\ge_{\mathrm{cp}}0.
\]
By the maximality property of the geometric mean with respect to the CP order, we obtain
\[
\Psi_a\circ\Theta\circ\Psi_{a^{-1}} \le_{\mathrm{cp}} \Theta,
\]
and thus,
\[
\Psi_a\circ\Theta \le_{\mathrm{cp}} \Theta\circ\Psi_a .
\]
Exchanging $a$ and $a^{-1}$ yields the reverse inequality, so
\[
\Theta\circ\Psi_a=\Psi_a\circ\Theta.
\]
Thus
\[
\Theta(axa^*)=a\Theta(x)a^* .
\]

For unitaries $u, v\in \cN_1\cap\cN_2$ and $z\in\C$,
we put $a=u+zv=u(1+zu^*v)$.
If $|z|<1$, then $\|zu^*v\|<1$,
and so $1+zu^*v$ is invertible.
Thanks to $\Theta(axa^*)=a\Theta(x)a^*$,
we can obtain
\[
\Theta(vxu^*)=v\Theta(x)u^*.
\]
For any $a, b\in \cN_1\cap\cN_2$,
using the standard fact that every element in a C$^*$-algebra is a linear combination of unitaries,
we conclude 
\[
\Theta(axb)=a\Theta(x)b.
\]

(2)
By Corollary~\ref{cor:unital-bound}, the geometric mean of two unital CP maps 
is subunital, hence $\Theta(1)\le1$.

Applying the bimodule property to $x=1$ gives
\[
a\Theta(1)=\Theta(a)=\Theta(1)a
\quad\text{for}\ a\in \cN_1\cap\cN_2,
\]
which shows
\[
\Theta(1)\in (\cN_1\cap\cN_2)' \cap \cM .
\]

\epf

We can similarly prove the following:

\bcor
Let $E_1, E_2 \colon \cM \to \cN$ be conditional expectations.
Then
\benu
\item $E_1\# E_2$ is an $\cN$-bimodule map, 
\item $E_1\# E_2$ is subunital and satisfies
\[
(E_1\# E_2)(1) \in \cN' \cap \cM.
\]
\eenu
\ecor

The notion of index for CP maps was introduced 
by Pimsner and Popa \cite{pp}, 
and further studied in \cite{izu}.
The Pomsner--Popa index $\mathrm{Ind}_{\mathrm{cp}}(\Phi)$ of $\Phi\in\mathrm{CP}(\cM, \cM)$ is defined by
\[
\mathrm{Ind}_{\mathrm{cp}}(\Phi)\coloneqq
\inf\{\lambda>0\mid
\lambda\Phi-\id_{\cM}\geq_{\mathrm{cp}} 0\}.
\]

\bprop\label{prop:pimsner-popa}
Let $\Phi_i\in\mathrm{CP}(\cM, \cM)$ for $i=1,2$. 
If there exist $\lambda_i>0$ such that
\[
\Phi_i \geq_{\mathrm{cp}} \lambda_i\ \mathrm{id}_{\cM}
\quad\text{for}\ i=1,2,
\]
then
\[
\Phi_1 \# \Phi_2 \geq_{\mathrm{cp}} \sqrt{\lambda_1 \lambda_2}\ \mathrm{id}_\cM.
\]
In particular, if $\mathrm{Ind}_{\mathrm{cp}}(\Phi_i)<\infty$ 
for $i=1,2$,
then
\[
\mathrm{Ind}_{\mathrm{cp}}(\Phi_1\#\Phi_2)\leq\sqrt{\mathrm{Ind}_{\mathrm{cp}}(\Phi_1)\,\mathrm{Ind}_{\mathrm{cp}}(\Phi_2)}.
\]
\eprop

\bpf
Since $\Phi_i \geq_{\mathrm{cp}} \lambda_i \mathrm{id}_{\cM}$ 
for $i=1,2$, by Theorem \ref{thm:main-geometric},
we have
\[
\Phi_1 \# \Phi_2 \geq_{\mathrm{cp}} (\lambda_1 \mathrm{id}_{\cM}) \# (\lambda_2 \mathrm{id}_{\cM}) = \sqrt{\lambda_1 \lambda_2} \, \mathrm{id}_{\cM}.
\]
Hence, the statement holds.
\epf

\bex\label{ex:ce-tensor}
Let $\cN_1 = \cN_2 = M_n(\mathbb{C})$,
and $\cM=\cN_1\otimes\cN_2$.
Let $\psi_1$, $\psi_2$ be states on $\cN_1$, $\cN_2$ 
given by 
\[
\psi_1(x) = \mathrm{Tr}(\rho x),
\quad\psi_2(y) = \mathrm{Tr}(\sigma y)
\] 
with diagonal invertible density matrices 
\[
\rho=\mathrm{diag}(\rho_1,\dots,\rho_n),
\qquad
\sigma=\mathrm{diag}(\sigma_1,\dots,\sigma_n),
\] respectively. 
Then we define conditional expectations $E_1$, $E_2$ from $\cM$ onto $\cN_1$, $\cN_2$ by
\[
E_1(x_1\otimes x_2)=x_1\otimes \psi_2(x_2)1,
\quad
E_2(x_1\otimes x_2)=\psi_1(x_1)1\otimes x_2,
\]
respectively.

We consider the geometric mean $\tilde{E}_1\# \tilde{E}_2$.
We now compute the corresponding Choi matrices.
\begin{align*}
C_{\tilde{E}_1}
&=
\sum_{i,j,k,l}(e_{ij}\otimes e_{kl})\otimes \tilde{E}_1(e_{ij}\otimes e_{kl})
\\
&=
\sum_{i,j,k,l}(e_{ij}\otimes e_{kl})\otimes (e_{ij}\otimes \psi_2(e_{kl})1)
\\
&=
(\sum_{i,j}e_{ij}\otimes e_{ij})\otimes
(\sum_ke_{kk}\otimes \sigma_k1)
\\
&=
C_{\mathrm{id}_{\cN_1}}\otimes\boldsymbol{\sigma},
\end{align*}
where
$\boldsymbol{\sigma}$ is the $n^2 \times n^2$ diagonal matrix defined by
\[
\boldsymbol{\sigma} = \bigoplus_{k=1}^n \sigma_k I_n
\]
so that its diagonal entries are $(\underbrace{\sigma_1, \dots, \sigma_1}_{n}, \dots, \underbrace{\sigma_n, \dots, \sigma_n}_{n})$.
Similarly, we have
\[
C_{\tilde{E}_2}=\boldsymbol{\rho}\otimes C_{\mathrm{id}_{\cN_2}}.
\]
Notice that 
\[
\ip{v, A^{-1}v}^{-1}=\max\{\lambda>0
\mid A-\lambda vv^*\geq 0 \}
\]
for an invertible positive semidefinite matrix $A$ and a nonzero vector $v$.
Note that
\begin{align*}
0&\leq C_{\tilde{E}_1}-C_{\lambda\id_{\cM}}
\\
&=
C_{\mathrm{id}_{\cN_1}}\otimes\boldsymbol{\sigma}-\lambda C_{\mathrm{id}_{\cN_1}}\otimes C_{\mathrm{id}_{\cN_2}}
\\
&=
C_{\mathrm{id}_{\cN_1}}\otimes(\boldsymbol{\sigma}-\lambda C_{\mathrm{id}_{\cN_2}}).
\end{align*}
Recall that $C_{\mathrm{id}}=vv^*$, where
\[
v=\sum_i e_i\otimes e_i.
\]
Hence
\[
\lambda_\sigma\coloneqq\mathrm{Ind}_{\mathrm{cp}}(\tilde{E}_1)^{-1}=(\sum_i\sigma_i^{-1})^{-1}
=
\ip{
v, \boldsymbol{\sigma}^{-1}v
}^{-1}.
\]
Therefore
\[
\tilde{E}_1\geq_{\mathrm{cp}}\lambda_\sigma\id_{\cM}.
\]
Similarly,
\[
\tilde{E}_2\geq_{\mathrm{cp}}\lambda_\rho\id_{\cM}.
\]
By Proposition \ref{prop:pimsner-popa}, we have
\[
\tilde{E}_1\# \tilde{E}_2\geq_{\mathrm{cp}}\sqrt{\lambda_\rho\lambda_\sigma}\ \mathrm{id}_{\cM}.
\]
Moreover, 
by using Proposition \ref{prop:choi-matrix} 
and {\cite[Theorem 13]{and}},
we have
\begin{align*}
C_{\tilde{E}_1\#\tilde{E}_2}
&=
C_{\tilde{E}_1}\#C_{\tilde{E}_2}
\\
&=
(C_{\mathrm{id}_{\cN_1}}\otimes\boldsymbol{\sigma})\#(\boldsymbol{\rho}\otimes C_{\mathrm{id}_{\cN_2}})
\\
&=
(C_{\mathrm{id}_{\cN_1}}\#\boldsymbol{\rho})\otimes(\boldsymbol{\sigma}\# C_{\mathrm{id}_{\cN_2}}).
\end{align*}
Since
\[
\boldsymbol{\rho}^{-1/2}C_{\mathrm{id}}\boldsymbol{\rho}^{-1/2}
=
(\boldsymbol{\rho}^{-1/2}v)(\boldsymbol{\rho}^{-1/2}v)^*,
\]
we have
\[
(\boldsymbol{\rho}^{-1/2}C_{\mathrm{id}_{\cN_1}}\boldsymbol{\rho}^{-1/2})^{1/2}
=\frac{(\boldsymbol{\rho}^{-1/2}v)(\boldsymbol{\rho}^{-1/2}v)^*}{\sqrt{\sum_{i=1}^n\rho_i^{-1}}}.
\]
Therefore
\[
C_{\mathrm{id}_{\cN_1}}\#\boldsymbol{\rho}
=
\boldsymbol{\rho}^{1/2}(\boldsymbol{\rho}^{-1/2}C_{\mathrm{id}_{\cN_1}}\boldsymbol{\rho}^{-1/2})^{1/2}\boldsymbol{\rho}^{1/2}
=
\frac{1}{\sqrt{\sum_{i=1}^n\rho_i^{-1}}}C_{\mathrm{id}_{\cN_1}}.
\]
Similarly,
\[
\boldsymbol{\sigma}\# C_{\mathrm{id}_{\cN_2}}
=
\frac{1}{\sqrt{\sum_{i=1}^n\sigma_i^{-1}}}C_{\mathrm{id}_{\cN_2}}.
\]
Then we obtain
\[
C_{\tilde{E}_1\#\tilde{E}_2}
=
\sqrt{\lambda_\rho\lambda_\sigma}\ (C_{\mathrm{id}_{\cN_1}}\otimes C_{\mathrm{id}_{\cN_2}})
=
\sqrt{\lambda_\rho\lambda_\sigma}\ C_{\mathrm{id}_{\cM}}.
\]
Therefore
\[
\tilde{E}_1\#\tilde{E}_2=\sqrt{\lambda_\rho\lambda_\sigma}\ \mathrm{id}_{\cM},
\]
which implies
\[
\mathrm{Ind}_{\mathrm{cp}}(\tilde{E}_1\#\tilde{E}_2)
=
(\sqrt{\lambda_\rho\lambda_\sigma})^{-1}
=
\sqrt{(\sum_{i=1}^n\rho_i^{-1})(\sum_{i=1}^n\sigma_i^{-1})}.
\]
\eex

\bex
Let $\cM= M_2(\C)\supset \cN_1=\C^2$ (diagonal),
and $E_1\colon\cM\to\cN_1$ be the canonical conditional expectation.
Let $\theta\in\R$ and set
\[
u=
\begin{bmatrix}
\cos\theta & -\sin\theta \\
\sin\theta & \cos\theta \\
\end{bmatrix},
\]
and $\cN_2=u\cN_1 u^*$.
Then we have the conditional expectation 
$E_2\colon\cM\to\cN_2$ by
\[
E_2(x)=u E_1(u^*xu)u^*
\quad\text{for}\ x\in\cM.
\]
To compute the geometric mean $\tilde{E}_1\#\tilde{E}_2$,
we consider the corresponding Choi matrices.
One can check that 
\[
C_{\tilde{E}_1}
=
\sum_{i}e_{ii}\otimes e_{ii}
=
\begin{bmatrix}
1 & 0 & 0 & 0 \\
0 & 0 & 0 & 0 \\
0 & 0 & 0 & 0 \\
0 & 0 & 0 & 1 \\
\end{bmatrix}
\] 
is the projection
onto $\cV_1=\mathrm{span}\{e_1\otimes e_1, e_2\otimes e_2\}$.
Since $u$ is real,
\[
C_{\tilde{E}_2}=(u\otimes u)C_{\tilde{E}_1}
(u^*\otimes u^*)
\] 
is the projection
onto $\cV_2=\mathrm{span}\{ue_1\otimes ue_1, ue_2\otimes ue_2\}$.
Note that
\[
\cV_1\cap\cV_2=
\begin{cases}
\mathrm{span}\{e_1\otimes e_1+e_2\otimes e_2\} & (\sin(2\theta)\ne 0), \\
\cV_1 & (\sin(2\theta)= 0).
\end{cases}
\]
Let $\Theta=\tilde{E}_1\#\tilde{E}_2$.
Since $C_{\tilde{E}_1}$ and $C_{\tilde{E}_2}$ are projections,
we obtain
\[
C_{\Theta} = C_{\tilde{E}_1} \# C_{\tilde{E}_2}=
C_{\tilde{E}_1} \wedge C_{\tilde{E}_2}.
\]
If $\sin(2\theta)\ne 0$, 
then 
\[
C_\Theta=\frac{1}{\|v\|^2}vv^*=\frac{1}{2}\sum_{i,j}e_{ij}\otimes e_{ij}=\frac{1}{2}C_{\id_\cM},
\]
where
$v=e_1\otimes e_1+e_2\otimes e_2$.
Therefore, by the Choi--Jamio\l{}kowski isomorphism,
we have
\[
\tilde{E}_1\# \tilde{E}_2
=
\begin{cases}
\frac{1}{2}\mathrm{id}_\cM & (\sin(2\theta)\ne 0),\\
\tilde{E}_1 & (\sin(2\theta)= 0). \\
\end{cases}
\]
\eex

\section{Connections for completely positive maps}

In this section, we introduce a notion of connections for CP maps.
Our goal is to extend the Kubo--Ando theory of operator connections to CP maps
in a way that is compatible with their natural order structure and composition.

We begin by recalling the notion of connections for bounded positive operators on an infinite-dimensional complex Hilbert space $\cH$ in the Kubo--Ando theory \cite{ka}.

\bdf
A binary operation 
$\sigma : \bB(\cH)^+\times \bB(\cH)^+\to \bB(\cH)^+$ is called an {\em operator connection}
if it satisfies
\benu
\item $A\leq C$ and $B\leq D$ imply
$A\sigma B\leq C \sigma D$,
\item
$C(A\sigma B)C\leq (CAC)\sigma (CBC)$,
\item
$A_n\downarrow A$ and $B_n\downarrow B$ 
imply $A_n\sigma B_n\downarrow A\sigma B$
in the strong operator topology.
\eenu

An operator connection $\sigma$
is called an {\em operator mean} if 
\begin{enumerate}
\item[(4)]
$I\sigma I=I$.
\end{enumerate}
\edf

\bdf
We denote by $\mathrm{OM}^+([0,\infty))$
the set of all operator monotone functions $f\colon[0,\infty)\to[0,\infty)$.
\edf

\bthm[Kubo--Ando]
For each operator connection $\sigma$,
there exists a unique 
$f\in \mathrm{OM}^+([0,\infty))$
such that 
\[
f(t)I=I\sigma (tI)
\]
for $t\geq 0$.
Moreover,
the map $\sigma\mapsto f$ is an affine order-isomorphism
such that 
$\sigma$ is an operator mean
if and only if $f(1)=1$.
In this case, $A\sigma A=A$
for all $A\in \bB(\cH)^+$.
\ethm

By L\"owner theory, 
for each $f\in \mathrm{OM}^+([0,\infty))$,
there exist unique constants $a, b\geq 0$ and
a unique finite positive measure $\mu$ on $(0,\infty)$
such that 
\begin{equation}\label{eq:monotone}
f(t)=a+bt+\int_{(0,\infty)}
\frac{t(1+\lambda)}{t+\lambda}\,d\mu(\lambda).
\end{equation}
Then the corresponding operator connection
$\sigma_f$ is given by
\[
A\sigma_f B=aA+bB+\int_{(0,\infty)}
\frac{1+\lambda}{\lambda}\Big[(\lambda A) : B\Big]\,d\mu(\lambda).
\]
Moreover, for $\alpha, \beta\in F_+(\cV)$,
we define the form connection by
\[
\alpha\sigma_f \beta(x,y)
\coloneqq a\alpha(x,y)+b\beta(x,y)+\int_{(0,\infty)}
\frac{1+\lambda}{\lambda}\Big[(\lambda \alpha) : \beta\Big](x,y)\,d\mu(\lambda).
\]
It is easy to see that $\alpha\sigma_f \beta\in F_+(\cV)$.

We refer to \cite{hk2} for the notion of connections for positive sesquilinear forms, 
formulated in terms of (possibly unbounded) quadratic forms. 
While their setting differs from ours, it provides a useful conceptual background for the development of connections for CP maps.

Let $f\in \mathrm{OM}^+([0,\infty))$.
For $\Phi, \Psi\in\mathrm{CP}(\cM, \cN)$,
we consider the corresponding positive sesquilinear forms
$s_\Phi, s_\Psi$ on $\cM\odot L^2(\cN)$.
Then we set 
\[
\gamma\coloneqq s_\Phi\sigma _fs_\Psi.
\]
Via separation and completion
of $\cM\odot L^2(\cN)$ with respect to $\gamma$,
we obtain the Hilbert space $\cH_\gamma$.
Then $\cH_\gamma$ is equipped with commuting left and right actions of $\cM$ and $\cN$, respectively.
For $\xi\in L^2(\cN)$, we have
\begin{align*}
&
\int_{(0,\infty)}
\frac{1+\lambda}{\lambda}\Big[(\lambda s_\Phi) : s_\Psi\Big](1\otimes \xi, 1\otimes \xi)
\,d\mu(\lambda)
\\
&\quad\leq
\int_{(0,1]}
\frac{1+\lambda}{\lambda}\lambda s_\Phi(1\otimes \xi, 1\otimes \xi)
\,d\mu(\lambda)
+
\int_{[1,\infty)}
\frac{1+\lambda}{\lambda}s_\Psi(1\otimes \xi, 1\otimes \xi)
\,d\mu(\lambda)
\\
&\quad\leq
2(\|\Phi\|+\|\Psi\|) \|\mu\| \|\xi\|^2
\end{align*}
Therefore,
\begin{align*}
\gamma(1\otimes \xi, 1\otimes \xi)
&=
as_\Phi(1\otimes \xi, 1\otimes \xi)+bs_\Psi(1\otimes \xi, 1\otimes \xi)\\
&\quad+\int_{(0,\infty)}
\frac{1+\lambda}{\lambda}[(\lambda s_\Phi) : s_\Psi](1\otimes \xi, 1\otimes \xi)\,d\mu(\lambda)
\\
&\leq
\Big\{a\|\Phi\|+b\|\Psi\|+2(\|\Phi\|+\|\Psi\|) \|\mu\|\Big\}
\|\xi\|^2.
\end{align*}
Hence, the map 
\[
V : L^2(\cN)\to \cH_\gamma,
\quad
V(\xi)=h(1\otimes\xi)
\]
extends to a bounded operator.
We define the CP map
\[
\Theta(x)=V^*xV
\]
for $x\in \cM$.
Moreover, 
one can check that $\gamma$ satisfies (\ref{eq:right-action}).
Hence, for $b\in\cN$, 
\[
\ip{
\Theta(x)Jb^*J\xi, \eta
}
=
\ip{
Jb^*J\Theta(x)\xi, \eta
}.
\]
Therefore, 
\[
\Theta(x)\in (\cN')'=\cN,
\]
which implies $\Theta\in\mathrm{CP}(\cM, \cN)$.
By (\ref{eq:left-action}),
we also have
\[
\ip{
\Theta(y^*x)\xi, \eta
}
=
(s_\Phi \sigma_f s_\Psi)(x\otimes \xi, y\otimes \eta).
\]

\bdf
For $f\in \mathrm{OM}^+([0,\infty))$
given in (\ref{eq:monotone}),
the {\em CP-map connection} $\sigma_f$
associated with $f$ is defined by
\[
\Phi\sigma_f\Psi\coloneqq\Theta,
\qquad
\Phi,\Psi\in\mathrm{CP}(\cM,\cN).
\]
\edf

Since $\gamma$ also satisfies (\ref{eq:left-action}), 
we obtain 
\[
s_{\Phi\sigma_f\Psi}=\gamma=s_\Phi\sigma_f s_\Psi.
\]
We remark that if $\Phi$, $\Psi$ are normal,
then so is $\Phi\sigma_f\Psi$.

\bprop
Let $f_1, f_2\in \mathrm{OM}^+([0,\infty))$.
Then
$\Phi\sigma_{f_1}\Psi\leq_{\mathrm{cp}}
\Phi\sigma_{f_2}\Psi$ for all $\Phi, \Psi\in\mathrm{CP}(\cM, \cN)$
if and only if $f_1(t)\leq f_2(t)$ for all $t\geq 0$.
\eprop

\bpf
We realize $s_\Phi, s_\Psi$ on the Hilbert space $\cH_\gamma$ with respect to $\gamma=s_\Phi+s_\Psi$
as
\[
s_\Phi(x\otimes \xi, y\otimes \eta)=\langle Ah(x\otimes \xi), h(y\otimes \eta)\rangle,
\]
and
\[
s_\Psi(x\otimes \xi, y\otimes \eta)=\langle Bh(x\otimes \xi), h(y\otimes \eta)\rangle.
\]
Then one checks that
\[
s_\Phi\sigma s_\Psi(\zeta, \zeta)
=
\langle (A\sigma B)h(\zeta),h(\zeta)\rangle.
\]
Hence the assertion follows from the operator case.
\epf

Using the corresponding properties of the parallel sum,
one can verify that every CP-map connection satisfies the following properties.

\bthm
Let $\sigma$ be a CP-map connection
and $f$ be the corresponding function.
\benu
\item $\Phi_1\leq_{\mathrm{cp}}\Phi_2$ and
$\Psi_1\leq_{\mathrm{cp}}\Psi_2$ imply
$\Phi_1\sigma\Psi_1\leq_{\mathrm{cp}}\Phi_2\sigma\Psi_2$,
\item For $\Xi\in\mathrm{CP}(\cN, \cP)$,
\[
\Xi\circ(\Phi\sigma_f\Psi)
\leq_{\mathrm{cp}}
(\Xi\circ\Phi)\sigma(\Xi\circ\Psi),
\]
\item For $\Xi\in\mathrm{CP}(\cL, \cM)$,
\[
(\Phi\sigma\Psi)\circ\Xi
\leq_{\mathrm{cp}}
(\Phi\circ\Xi)\sigma(\Psi\circ\Xi),
\]
\item $\Phi_n\downarrow\Phi$ and $\Psi_n\downarrow\Psi$ imply
$\Phi_n\sigma\Psi_n\downarrow\Phi\sigma\Psi$,
where $\Phi_n \downarrow \Phi$ means that 
$\Phi_1 \geq_{\mathrm{cp}} \Phi_2 \geq_{\mathrm{cp}} \cdots$ 
and $\Phi_n(x) \to \Phi(x)$ 
in the ultraweak operator topology for each $x \in \cM^+$.
\item If $f(1)=1$, then $\Phi\sigma\Phi=\Phi$.
\eenu
\ethm

\bcor
Let $\sigma$ be a CP-map connection
and $f$ be the corresponding function.
\benu
\item For $\alpha\in\mathrm{Aut}(\cN)$ and
$\beta\in\mathrm{Aut}(\cM)$,
\[
\alpha\circ(\Phi\sigma\Psi)\circ\beta
=
(\alpha\circ\Phi\circ\beta)\sigma(\alpha\circ\Psi\circ\beta),
\]
\item
\[
(\Phi_1\sigma\Psi_1)+(\Phi_2\sigma\Psi_2)
\leq_{\mathrm{cp}}
(\Phi_1+\Phi_2)\sigma(\Psi_1+\Psi_2).
\]
\eenu
\ecor

\bex
\benu
\item For $0\leq\alpha\leq 1$,
let $\#_\alpha$ denote the {\em $\alpha$-power mean},
which is the CP-map mean corresponding to the operator monotone function $t^\alpha$.
Note that $\Phi\#_0\Psi=\Phi$, $\Phi\#_1\Psi=\Psi$
and $\Phi\#_{1/2} \Psi=\Phi\#\Psi$.
\item The CP-map mean corresponding to the operator monotone function $(t-1)/\log t$ is called the {\em logarithmic mean} and denoted by $\Phi\lambda\Psi$.
\eenu
\eex

We give another proof of Theorem \ref{thm:am-gm-hm}.

\bcor
For $\Phi, \Psi\in\mathrm{CP}(\cM, \cN)$,
\[
\Phi\, !\, \Psi
\leq_{\mathrm{cp}}
\Phi\, \#\, \Psi
\leq_{\mathrm{cp}}
\Phi\, \lambda\, \Psi
\leq_{\mathrm{cp}}
\Phi\, \nabla\, \Psi.
\]
\ecor

\bpf
It immediately follows from the inequality
\[
\frac{2t}{1+t}\leq t^{1/2}\leq\frac{t-1}{\log t}\leq\frac{1+t}{2}
\]
for $t>0$.
\epf

\bthm
Let $\Phi_k, \Psi_k\in\mathrm{CP}(\cM_k,\cN_k)$
for $k=1,2$.
Then
\[
(\Phi_1\otimes \Phi_2)\#_\alpha(\Psi_1\otimes \Psi_2)
=
(\Phi_1\#_\alpha\Psi_1)\otimes(\Phi_2\#_\alpha\Psi_2).
\]
\ethm

\bpf
Set 
\[
\gamma_k=s_{\Phi_k}+s_{\Psi_k}
\]
for $k=1,2$.
By the standard construction via separation and completion, 
we have Hilbert spaces $\cH_{\gamma_k}$
and the canonical maps $h_k\colon \cM_k\odot L^2(\cN_k)\to \cH_{\gamma_k}$.
Then we have bounded operators $A_k, B_k$
on $\cH_{\gamma_k}$ such that
\[
s_{\Phi_k}(x_k\otimes\xi_k, y_k\otimes \eta_k)
=
\ip{
A_kh_k(x_k\otimes\xi_k), h_k(y_k\otimes \eta_k)
}_{\gamma_k},
\]
and
\[
s_{\Psi_k}(x_k\otimes\xi_k, y_k\otimes \eta_k)
=
\ip{
B_kh_k(x_k\otimes\xi_k), h_k(y_k\otimes \eta_k)
}_{\gamma_k}.
\]

Let $\Phi=\Phi_1\otimes \Phi_2$ and $\Psi=\Psi_1\otimes \Psi_2$.
We define a sesquilinear form $\gamma$ on
$(\cM_1 \odot L^2(\cN_1)) \odot (\cM_2 \odot L^2(\cN_2))$
by
\[
\gamma(\zeta_1 \otimes \zeta_2, \zeta'_1 \otimes \zeta'_2)
=
\gamma_1(\zeta_1, \zeta'_1)\gamma_2(\zeta_2, \zeta'_2)
\]
for $\zeta_1, \zeta_1'\in\cM_1\odot L^2(\cN_1)$
and $\zeta_2, \zeta_2'\in\cM_2\odot L^2(\cN_2)$.
Via separation and completion,
we have a Hilbert space $\cH_{\gamma}$.
Here we implicitly use the canonical identification
\[
(\cM_1\odot L^2(\cN_1))\odot (\cM_2\odot L^2(\cN_2))=(\cM_1\odot \cM_2)\odot (L^2(\cN_1)\odot L^2(\cN_2))
\]
and then $\cH_\gamma$ can be identified with
$\cH_{\gamma_1}\otimes \cH_{\gamma_2}$.

We set $h=h_1\otimes h_2$ and $A=A_1\otimes A_2$, $B=B_1\otimes B_2$.
Then 
it is straightforward to verify that
$s_\Phi \sim (h,A)$ and $s_\Psi \sim (h,B)$.
Moreover,
since $A_k$ and $B_k$ commute for each $k$,
$A$ and $B$ also commute.
Let 
\[x=x_1\otimes x_2,\quad y=y_1\otimes y_2\in\cM_1\otimes\cM_2
\]
and 
\[
\xi=\xi_1\otimes\xi_2,\quad \eta=\eta_1\otimes\eta_2\in L^2(\cN_1)\otimes L^2(\cN_2).
\]
By Pusz--Woronowicz theory,
\begin{align*}
&
(s_\Phi\#_\alpha s_\Psi)(x\otimes \xi, y\otimes \eta)
\\
&=
\ip{
A^{1-\alpha} B^\alpha h(x\otimes \xi), h(y\otimes \eta)
}_\gamma
\\
&=
\ip{
A_1^{1-\alpha} B_1^\alpha h_1(x_1\otimes \xi_1), h_1(y_1\otimes \eta_1)
}_{\gamma_1}
\ip{
A_2^{1-\alpha} B_2^\alpha h_2(x_2\otimes \xi_2), h_2(y_2\otimes \eta_2)
}_{\gamma_2}
\\
&=
(s_{\Phi_1}\#_\alpha s_{\Psi_1})(x_1\otimes \xi_1, y_1\otimes \eta_1)
(s_{\Phi_2}\#_\alpha s_{\Psi_2})(x_2\otimes \xi_2, y_2\otimes \eta_2)
\\
&=
(s_{\Phi_1\#_\alpha \Psi_1}\otimes s_{\Phi_2\#_\alpha\Psi_2})(x\otimes \xi, y\otimes \eta).
\end{align*}
Therefore,
\[
s_{\Phi\#_\alpha \Psi}=s_{\Phi_1\#_\alpha \Psi_1}\otimes s_{\Phi_2\#_\alpha\Psi_2}.
\]
Since the associated sesquilinear forms coincide, we conclude that
\[
\Phi\#_\alpha \Psi
=
(\Phi_1\#_\alpha \Psi_1)\otimes(\Phi_2\#_\alpha \Psi_2).
\]
\epf

\bdf
Let $\sigma$ be a CP-map mean
and $f$ be the corresponding function.
\benu
\item The CP-map mean with the representing function $tf(t^{-1})$ is called the {\em transpose} of $\sigma$
and denoted by $\sigma'$.
If $\sigma=\sigma'$, then $\sigma$ is said to be {\em symmetric}.
\item The CP-map mean with the representing function $f(t^{-1})^{-1}$ is called the {\em adjoint} of $\sigma$
and denoted by $\sigma^*$.
\item The CP-map mean with the representing function $t/f(t)$ is called the {\em dual} of $\sigma$
and denoted by $\sigma^\perp$.
\eenu
\edf

The following propositions follow from the definitions
and the corresponding results for operators.
These correspond exactly to the classical notions in Kubo--Ando theory.

\bprop
Let $\sigma$ be a CP-map mean
and $f$ be the corresponding function.
\benu
\item $\Phi\sigma\Psi=\Psi\sigma'\Phi$,
\item $(\sigma')'=\sigma$, $(\sigma^*)^*=\sigma$, and $(\sigma^\perp)^\perp=\sigma$.
\item $\sigma^\perp=(\sigma')^*=(\sigma^*)'$,
$\sigma'=(\sigma^*)^\perp=(\sigma^\perp)^*$,
and
$\sigma^*=(\sigma')^\perp=(\sigma^\perp)'$.
\eenu
\eprop

\bprop
If $\sigma$ is a symmetric mean,
then $!\leq_{\mathrm{cp}}\sigma\leq_{\mathrm{cp}}\nabla$.
\eprop

\bprop
For every CP-map mean $\sigma$, the following hold:
\benu
\item $(\Phi\sigma\Psi)+(\Psi\sigma\Phi)\leq_{\mathrm{cp}}\Phi+\Psi$,
\item $(\Phi\sigma\Psi)\, :\, (\Psi\sigma\Phi)\geq_{\mathrm{cp}}\Phi\, :\, \Psi$,
\item $(\Phi\sigma\Psi)\#(\Phi\sigma^\perp\Psi)=(\Phi\#\Psi)$,
\item $(\Phi\sigma\Psi)+(\Phi\sigma^\perp\Psi)\leq_{\mathrm{cp}}\Phi+\Psi$,
\item $(\Phi\sigma\Psi)\, :\, (\Phi\sigma^\perp\Psi)\geq_{\mathrm{cp}} (\Phi\, :\, \Psi)$.
\eenu
\eprop

\brem
In the present paper, we pay particular attention to the geometric mean $\#$ corresponding to $f(t) = \sqrt{t}$. The primary reason for this focus is its unparalleled symmetry. While the arithmetic mean $\nabla$ and the harmonic mean $!$ are symmetric (i.e., $\nabla = \nabla'$ and $! =\ !'$\,), they are mutually adjoint ($\nabla^* =\ !$\,) and dual ($\nabla^\perp =\ !$\,). Conversely, the geometric mean is identically self-adjoint ($\# = \#^*$) and self-dual ($\# = \#^\perp$). Thus, $\#$ serves as the unique invariant fixed point under all the fundamental transformations established in the propositions above.
\erem


\section{Lebesgue decomposition of completely positive maps}


In this section, we establish a Lebesgue-type decomposition for CP maps. 
Our approach is inspired by the classical theory for positive operators developed by Ando \cite{and2} 
and the form-theoretic framework of Simon \cite{Simon_1977}, 
and is further informed by subsequent refinements due to Kosaki \cite{kos_1985, kos_2018, kos_2023} 
and, more recently, Aibara and Ueda \cite{Aibara-Ueda}.

Related Radon--Nikodym type results for CP maps originate from the work of Arveson \cite{arv} 
and have been further studied in various contexts, including quantum operations and covariant CP maps 
(see, e.g., \cite{raginsky, joita}). 

The aim of this section is to establish such a decomposition within our framework, 
based on the parallel sum operation introduced in the previous sections.

More precisely, we introduce notions of absolute continuity and singularity for CP maps, 
and show that every CP map admits a decomposition into an absolutely continuous part 
and a singular part with respect to a fixed reference map. 

First, we introduce the notions of absolute continuity and singularity for CP maps based on the parallel sum. In what follows, the limits of CP maps are taken in the point-ultraweak operator topology.

\bdf
Let $\Phi, \Psi\in\mathrm{CP}(\cM, \cN)$.
\benu
\item $\Psi$ is said to be $\Phi$-{\em absolutely continuous}, denoted by $\Psi \ll \Phi$, if there exists a sequence of CP maps $\{\Psi_n\}$ and positive real numbers $\{\alpha_n\}$ such that 
\[
\Psi_n \uparrow \Psi \quad \text{and} \quad \Psi_n \leq_{\mathrm{cp}} \alpha_n \Phi
\quad\text{for all}\ n.
\]
\item $\Psi$ is said to be $\Phi$-{\em singular}, denoted by $\Phi \perp \Psi$, if for any $\Theta\in\mathrm{CP}(\cM, \cN)$,
\[
\Theta \leq_{\mathrm{cp}} \Phi \quad \text{and} \quad \Theta \leq_{\mathrm{cp}} \Psi \implies \Theta = 0. 
\]
\eenu
\edf

To prove the Lebesgue decomposition constructively, we use the parallel sum to obtain an algebraic characterization of these order-theoretic notions.

Let $\Phi, \Psi\in\mathrm{CP}(\cM, \cN)$.
Set $\Psi_n \coloneqq (n\Phi \, : \, \Psi)$. By Proposition~\ref{prop:parallel-basic} and Theorem~\ref{thm:main-parallel}, 
we have $\Psi_n \leq \Psi$ and the sequence $\{\Psi_n\}$ is monotonically increasing. 
Hence there exists the point-ultraweak limit 
\[
[\Phi]\Psi \coloneqq \lim_{n \to \infty} (n\Phi : \Psi)\in\mathrm{CP}(\cM, \cN)
\]
such that $[\Phi]\Psi\leq_\mathrm{cp}\Psi$.
Since $\Psi_n \leq n\Phi$ by Proposition~\ref{prop:parallel-basic},
setting $\alpha_n = n$, the sequence $\{\Psi_n\}$ satisfies the condition for $[\Phi]\Psi \ll \Phi$.

This provides a canonical choice of the absolutely continuous component.
For convenience, we also write
\[
\Psi_{\mathrm{ac}} := [\Phi]\Psi.
\]

\blem\label{lem:parallel_characterization}
Let $\Phi, \Psi\in\mathrm{CP}(\cM, \cN)$.
\benu
\item $\Phi \perp \Psi$ if and only if $\Phi \, : \, \Psi = 0$.
\item $\Psi \ll \Phi$ if and only if $\Psi =[\Phi]\Psi$.
\eenu
\elem

\bpf
(1) Suppose $\Theta \leq_{\mathrm{cp}} \Phi$ 
and $\Theta \leq_{\mathrm{cp}} \Psi$. 
By the monotonicity of the parallel sum (Theorem \ref{thm:main-parallel}), 
\[
\frac{1}{2}\Theta = \Theta \, : \, \Theta \leq_{\mathrm{cp}} \Phi \, : \, \Psi. 
\]
Hence, if $\Phi \, : \, \Psi = 0$, then $\Theta = 0$, 
which means $\Phi \perp \Psi$. 
Conversely, by Proposition \ref{prop:parallel-basic},
we have
\[
\Phi \, : \, \Psi \leq_{\mathrm{cp}} \Phi \quad\text{and} \quad \Phi \, : \, \Psi \leq_{\mathrm{cp}} \Psi. 
\] 
If $\Phi \perp \Psi$, then
it immediately follows that $\Phi \, : \, \Psi = 0$.

(2) 
If 
\[
\Psi = \lim_{n \to \infty} (n\Phi : \Psi),
\]
then $\Psi \ll \Phi$.

Conversely, assume that $\Psi \ll \Phi$. 
Then
there exists a sequence of CP maps $\{\Psi_m\}$ and positive real numbers $\{\alpha_m\}$ such that 
\[
\Psi_m \uparrow \Psi \quad \text{and} \quad \Psi_m \leq_{\mathrm{cp}} \alpha_m \Phi
\quad\text{for all}\ m.
\]

Now we claim that
\[
\Psi_m-(n\Phi : \Psi_m) \leq_{\mathrm{cp}} \frac{\alpha_m}{n+\alpha_m}\Psi. 
\]
Consider a sesquilinear form 
$\gamma=s_\Phi+s_\Psi$.
Then there exist compatible representations
$s_\Phi\sim(h, A)$ and $s_\Psi\sim(h,B)$ on the Hilbert space $\cH_\gamma$.
Then, since $\Psi_m\leq_{\mathrm{cp}}\Psi$,
there exist bounded operators $B_m$
such that 
\[
s_{\Psi_m}(x\otimes\xi, y\otimes\eta)
=
\ip{
B_mh(x\otimes\xi), h(y\otimes\eta)
}.
\]
Hence we obtain $B_m\uparrow B$ and $B_m\leq\alpha_m A$, i.e., $B$ is $A$-absolutely continuous.
By the proof of \cite[Lemma 1]{and2},
we have
\[
B_m-(nA : B_m)\leq\frac{\alpha_m}{n+\alpha_m}B.
\]
Therefore, we conclude our claim.

Taking the limit as $n \to \infty$, we obtain 
\[
\Psi_m=\lim_{n\to\infty}(n\Phi : \Psi_m),
\]
which means $\Psi_m\ll\Phi$ by the above.
Since
\[
\Psi_m=
\lim_{n\to\infty}(n\Phi : \Psi_m)
\leq_{\mathrm{cp}}
\lim_{n\to\infty}(n\Phi : \Psi)
\leq_{\mathrm{cp}}\Psi,
\]
taking the limit as $m \to \infty$,
we have 
\[
\Psi=\lim_{n\to\infty}(n\Phi : \Psi)=[\Phi]\Psi.
\]
\epf

\bthm\label{thm:Lebesgue Decomposition}
Let $\Phi, \Psi\in\mathrm{CP}(\cM, \cN)$. Then there exists a decomposition
\[
\Psi = \Psi_{\mathrm{ac}} + \Psi_{\mathrm{s}}
\]
where $\Psi_{\mathrm{ac}}, \Psi_{\mathrm{s}}\in\mathrm{CP}(\cM, \cN)$ satisfy $\Psi_{\mathrm{ac}}\ll \Phi$ and $\Psi_{\mathrm{s}} \perp \Phi$.
This decomposition is not unique in general. 
However, there exists a canonical largest choice for the absolutely continuous part.
Specifically, the maximum among all $\Theta\in\mathrm{CP}(\cM, \cN)$ such that $\Theta \leq_{\mathrm{cp}} \Psi$ and $\Theta \ll \Phi$ 
is given by
\[
\Psi_{\mathrm{ac}} = \lim_{n \to \infty} (n\Phi \, : \,\Psi). 
\]
Furthermore, this Lebesgue decomposition is unique if and only if there exists a constant $\alpha > 0$ such that $\Psi_{\mathrm{ac}} \leq_{\mathrm{cp}} \alpha \Phi$.
\ethm

\bpf
Recall that $\Psi_{\mathrm{ac}} \ll \Phi$.
We then define the remaining part as 
\[
\Psi_{\mathrm{s}} \coloneqq \Psi - \Psi_{\mathrm{ac}}\in \mathrm{CP}(\cM, \cN).
\]
Consider a sesquilinear form 
$\gamma=s_\Phi+s_\Psi$.
Then there exist compatible representations
$s_\Phi\sim(h, A)$ and $s_\Psi\sim(h,B)$ on the Hilbert space $\cH_\gamma$.
Note that
\[
\ip{
(n\Phi\, :\, \Psi)(y^*x)\xi, \eta
}
=
\ip{
(nA : B)h(x\otimes \xi), h(y\otimes \eta)
}_\gamma.
\]
Hence
\[
\ip{
\Psi_{\mathrm{ac}}(y^*x)\xi, \eta
}
=
\ip{
[A]Bh(x\otimes \xi), h(y\otimes \eta)
}_\gamma,
\]
where 
\[
[A]B=\lim_{n\to\infty}(nA):B.
\]
Therefore
\[
\ip{
\Psi_{\mathrm{s}}(y^*x)\xi, \eta
}
=
\ip{
(B-[A]B)h(x\otimes \xi), h(y\otimes \eta)
}_\gamma.
\]
By \cite{and2},
the corresponding operator $B-[A]B$ is $A$-singular.
If $\Theta\in\mathrm{CP}(\cM, \cN)$ such that $\Theta \leq_{\mathrm{cp}}\Phi$ and $\Theta \leq_{\mathrm{cp}} \Psi_{\mathrm{s}}$,
then there exists a positive operator $D$
such that
\[
\ip{
\Theta(y^*x)\xi, \eta
}
=
\ip{
Dh(x\otimes \xi), h(y\otimes \eta)
}_\gamma.
\]
Since $D\leq A$, $D\leq B-[A]B$,
we have $D=0$, which means $\Theta=0$.
Therefore $\Psi_{\mathrm{s}}$ is $\Phi$-singular.

Let $\Theta\in\mathrm{CP}(\cM, \cN)$ such that $\Theta \leq_{\mathrm{cp}} \Psi$ and $\Theta \ll \Phi$.
By Lemma \ref{lem:parallel_characterization},
we have
\[
\Theta=\lim_{n\to\infty}(n\Phi : \Theta)
\leq_{\mathrm{cp}}\lim_{n\to\infty}(n\Phi : \Psi)
=\Psi_{\mathrm{ac}}.
\]

Finally, we consider the condition for uniqueness. 
Suppose that there is another Lebesgue decomposition $\Psi = \Psi_1 + \Psi_2$ where $\Psi_1 \ll \Phi$ and $\Psi_2 \perp \Phi$. 
Passing to the Hilbert space $\cH_\gamma$, this corresponds to a decomposition of positive operators $B = B_1 + B_2$ such that $B_1 \ll A$ and $B_2 \perp A$. By \cite[Theorem 6]{and2}, 
this decomposition is unique (meaning $B_1 = [A]B$) if and only if $[A]B \leq \alpha A$ for some $\alpha\geq 0$. 
Translating this back to CP maps via the affine order-isomorphism, $[A]B \leq \alpha A$ holds if and only if $\Psi_{\mathrm{ac}} \leq_{\mathrm{cp}} \alpha \Phi$.
Thus, the Lebesgue decomposition for CP maps is unique if and only if $\Psi_{\mathrm{ac}} \leq_{\mathrm{cp}} \alpha \Phi$ for some $\alpha > 0$.
\epf

\bex[Lebesgue decomposition of normal positive functionals]
\label{ex:normal-states}
Let $(\cM,\cH,J,\cP)$ be a standard form of a von Neumann algebra
$\cM$, and let $\varphi,\psi\in\cM_*^+$.  Set
$\omega=\varphi+\psi$ and $r=s(\omega)$.  Both $\varphi$ and $\psi$
are supported by $r$, that is,
\[
 \varphi(x)=\varphi(rxr),
 \qquad
 \psi(x)=\psi(rxr),
 \qquad x\in\cM.
\]
If $\omega=0$, then $\varphi=\psi=0$ and there is nothing to prove,
so we assume that $\omega\neq0$.
Compression and extension by $x\mapsto rxr$ give an order
isomorphism between $(r\cM r)_*^+$ and the normal positive
functionals on $\cM$ supported by $r$.  This order isomorphism
preserves domination, absolute continuity, singularity, and the
maximal absolutely continuous component.  We may therefore pass to
the standard form of $r\cM r$ and relabel the corner as $\cM$.
After this reduction, $\omega$ is faithful.

Let $\xi_\omega\in\cP$ be the unique implementing vector for
$\omega$.  Since $\omega$ is faithful, $\xi_\omega$ is cyclic and
separating for $\cM$.  As $\varphi,\psi\leq\omega$, there exist unique
positive contractions $h'_\varphi,h'_\psi\in\cM'$ such that $h'_\varphi+h'_\psi=1$ and
\[
 \varphi(x)=\ip{h'_\varphi x\xi_\omega,\xi_\omega},
 \qquad
 \psi(x)=\ip{h'_\psi x\xi_\omega,\xi_\omega}.
 \]
Write
$e'_\varphi\coloneqq s(h'_\varphi)$
and
$e'_\psi\coloneqq s(h'_\psi)$.
Since $h'_\psi=1-h'_\varphi$, the operators $h'_\varphi$ and
$h'_\psi$ commute.  For $n\geq1$, define
\[
 f_n(t)\coloneqq\frac{nt(1-t)}{nt+1-t},
 \qquad 0\leq t\leq1.
\]
Then
$0\leq f_n\leq1$ and
\[
 f_n(t)\to
 \chi_{(0,1]}(t)(1-t).
\]
Moreover, 
\[
 nh'_\varphi:h'_\psi=f_n(h'_\varphi).
\]
By bounded Borel functional calculus,
we have
\begin{align*}
\lim_{n\to\infty}
       (nh'_\varphi:h'_\psi)
 &=\chi_{(0,1]}(h'_\varphi)(1-h'_\varphi)\\
 &=e'_\varphi h'_\psi
 =h'_\psi e'_\varphi.
\end{align*}
Consequently, 
\[
 \psi_{\mathrm{ac}}(x)=\ip{e'_\varphi h'_\psi x\xi_\omega,\xi_\omega},
 \qquad
 \psi_{\mathrm{s}}(x)
 =\ip{(1-e'_\varphi)h'_\psi x\xi_\omega,\xi_\omega}.
\]

We next compare this formula with Kosaki's construction
\cite{kos_1985}.  Put
\[
 a=[D\varphi:D\omega]_{-i/2}\in\cM,
\]
let $p$ be the projection onto $\ker a$, and set
$p'=JpJ\in\cM'$.  One has
\[
 \varphi(x)=\ip{Ja^*aJx\xi_\omega,\xi_\omega},
 \qquad x\in\cM,
\]
and hence
\[
 h'_\varphi=Ja^*aJ,
 \qquad
 1-p'=s(h'_\varphi)=e'_\varphi.
\]
Since $0\leq a^*a\leq 1$ and
$s(a^*a)=1-p$, one has the inequality
\[
 0\leq a^*a\leq 1-p.
\]
It follows that Kosaki's functional
\begin{align*}
 \widetilde\psi(x)
 &=\ip{(1-p')x\xi_\omega,\xi_\omega}-\varphi(x)\\
 &=\ip{J\bigl((1-p)-a^*a\bigr)Jx\xi_\omega,
              \xi_\omega}
\end{align*}
is positive.  Using $1-p'=e'_\varphi$ and
$h'_\psi=1-h'_\varphi$, we obtain
\begin{align*}
 \widetilde\psi(x)
 &=\ip{(e'_\varphi-h'_\varphi)x\xi_\omega,
             \xi_\omega}\\
 &=\ip{e'_\varphi h'_\psi x\xi_\omega,
             \xi_\omega}
 =\psi_{\mathrm{ac}}(x).
\end{align*}
Thus the parallel-sum construction agrees with Kosaki's maximal
$\varphi$-absolutely continuous part.  In addition,
\begin{align*}
 \ker a=\{0\}
 &\iff p'=0
 \iff e'_\varphi=1\\
 &\iff e'_\psi\leq e'_\varphi
 \iff \psi\ll\varphi.
\end{align*}
For the only non-immediate implication in the middle, if
$e'_\psi\leq e'_\varphi$, then
$\ker h'_\varphi\subseteq\ker h'_\psi$; since
$h'_\varphi+h'_\psi=1$, this forces
$\ker h'_\varphi=\{0\}$.

Returning to the original algebra, let
\[
 \varphi_r=\varphi|_{r\cM r},
 \qquad
 \psi_r=\psi|_{r\cM r}.
\]
Suppose that
\[
 \psi_r=(\psi_r)_{\mathrm{ac}}+(\psi_r)_{\mathrm{s}}
\]
is the Lebesgue decomposition obtained in the corner.  Define
normal positive functionals on $\cM$ by
\[
 \psi_{\mathrm{ac}}(x)
 =
 (\psi_r)_{\mathrm{ac}}(rxr),
 \qquad
 \psi_{\mathrm{s}}(x)
 =
 (\psi_r)_{\mathrm{s}}(rxr)
 \qquad x\in\cM.
\]
Since $\psi(x)=\psi_r(rxr)$ for every $x\in\cM$, we have
\[
 \psi_{\mathrm{s}}(x)
 =
 \psi(x)-\psi_{\mathrm{ac}}(x).
\]
By the preservation properties established above,
\[
 \psi_{\mathrm{ac}}\ll\varphi,
 \qquad
 \psi_{\mathrm{s}}\perp\varphi.
\]
For completeness, we verify these assertions directly.

Since $(\psi_r)_{\mathrm{ac}}\ll\varphi_r$, there exist
$\theta_n\in(r\cM r)_*^+$ and $\alpha_n>0$ such that
$\theta_n\uparrow(\psi_r)_{\mathrm{ac}}$
and
$\theta_n\leq\alpha_n\varphi_r$.
Define $\psi_n\in\cM_*^+$ by
\[
 \psi_n(x)\coloneqq\theta_n(rxr),
 \qquad x\in\cM.
\]
Then, for every $x\in\cM_+$,
\[
 \psi_n(x)
 =
 \theta_n(rxr)
 \uparrow
 (\psi_r)_{\mathrm{ac}}(rxr)
 =
 \psi_{\mathrm{ac}}(x),
\]
and
\[
 \psi_n(x)
 =
 \theta_n(rxr)
 \leq
 \alpha_n\varphi_r(rxr)
 =
 \alpha_n\varphi(x).
\]
Thus
\[
 \psi_n\uparrow\psi_{\mathrm{ac}},
 \qquad
 \psi_n\leq\alpha_n\varphi,
\]
and hence $\psi_{\mathrm{ac}}\ll\varphi$.

To prove singularity, suppose that $\theta\in\cM_*^+$ satisfies
$\theta\leq\varphi$
and
$\theta\leq\psi_{\mathrm{s}}$.
Let
\[
  \theta_r=\theta|_{r\cM r}.
\]
Then
$\theta_r\leq\varphi_r$
and
$\theta_r\leq(\psi_r)_{\mathrm{s}}$.
Since $(\psi_r)_{\mathrm{s}}\perp\varphi_r$ in the corner, it follows
that $ \theta_r=0$.  Moreover, since $ \theta\leq\varphi$ and
$s(\varphi)\leq r$, we have
\[
 s( \theta)\leq s(\varphi)\leq r.
\]
Therefore,
\[
  \theta(x)= \theta(rxr)= \theta_r(rxr)=0,
 \qquad x\in\cM,
\]
so $ \theta=0$.  This proves that
$\psi_{\mathrm{s}}\perp\varphi$.

Finally, suppose that $\theta\in\cM_*^+$ satisfies
$\theta\leq\psi$ and $\theta\ll\varphi$.
Since $\theta\leq\psi$ and $s(\psi)\leq r$, we have
$s(\theta)\leq r$.  Hence
\[
 \theta(x)=\theta_r(rxr),
 \qquad
 \theta_r\coloneqq\theta|_{r\cM r}.
\]
Moreover,
$\theta_r\leq\psi_r$ and $\theta_r\ll\varphi_r$.
By the maximality of $(\psi_r)_{\mathrm{ac}}$ in the corner,
\[
 \theta_r\leq(\psi_r)_{\mathrm{ac}}.
\]
Consequently, for every $x\in\cM_+$,
\[
 \theta(x)
 =
 \theta_r(rxr)
 \leq
 (\psi_r)_{\mathrm{ac}}(rxr)
 =
 \psi_{\mathrm{ac}}(x).
\]
Thus $\theta\leq\psi_{\mathrm{ac}}$, and hence
$\psi_{\mathrm{ac}}$ is the maximal $\varphi$-absolutely continuous
component of $\psi$ on the original algebra $\cM$.  No additional
component can occur outside the corner, since $\psi$ is supported
by $r$.
\eex

The preceding example shows that the Lebesgue decomposition is
encoded by the support projections of Radon--Nikodym derivatives.
The same description holds for arbitrary CP maps.

\bthm
Let $\Phi,\Psi\in\mathrm{CP}(\cM,\cN)$ and set
$\Gamma=\Phi+\Psi$.  Let $(\pi,V,\cH_\gamma)$ be a minimal
Stinespring representation of $\Gamma$, and let
$A',B'\in\pi(\cM)'$ be the Arveson Radon--Nikodym derivatives of
$\Phi$ and $\Psi$, respectively.  Thus
\[
 \Phi(x)=V^*A'\pi(x)V,
 \qquad
 \Psi(x)=V^*B'\pi(x)V,
  \qquad
  A'+B'=I.
\]
Let
\[
 P'_\Phi=s(A'),
 \qquad
 P'_\Psi=s(B'),
\]
and put
\[
 E_\Phi=P'_\Phi\cH_\gamma,
 \qquad
 E_\Psi=P'_\Psi\cH_\gamma.
\]
Then
\begin{align*}
 \Psi_{\mathrm{ac}}(x)
 &=V^*P'_\Phi B'P'_\Phi\pi(x)V,\\
 \Psi_{\mathrm{s}}(x)
 &=V^*(I-P'_\Phi)B'(I-P'_\Phi)\pi(x)V.
\end{align*}
Moreover,
\benu
\item
$\Psi\ll\Phi$ if and only if
$E_\Psi\subseteq E_\Phi$, equivalently
$P'_\Psi\leq P'_\Phi$;
\item
$\Psi\perp\Phi$ if and only if
$E_\Psi\perp E_\Phi$, equivalently
$P'_\Psi P'_\Phi=0$.
\eenu
\ethm

\bpf
Note that for every $n\geq1$,
\[
 (n\Phi:\Psi)(x)
 =V^*(nA':B')\pi(x)V.
\]
Since $A'+B'=I$, the operators $A'$ and $B'$ commute.  Hence
\[
 nA':B'
 =nA'B'(nA'+B')^{-1}\to
   P'_\Phi B'
 =P'_\Phi B'P'_\Phi.
\]
The strong limit follows from the functional calculus exactly
as in Example~\ref{ex:normal-states}.  Therefore
\[
 \Psi_{\mathrm{ac}}(x)
 =V^*P'_\Phi B'P'_\Phi\pi(x)V.
\]
Since $P'_\Phi$ commutes with $B'$,
we have
\begin{align*}
 B'-P'_\Phi B'P'_\Phi
 &=(I-P'_\Phi)B'(I-P'_\Phi),
\end{align*}
which proves the asserted formula for $\Psi_{\mathrm{s}}$.

We give a direct verification of the absolute-continuity criterion.
Suppose first that $P'_\Psi\leq P'_\Phi$.  For $k\geq1$, set
\[
 Q_k=\chi_{[1/k,1]}(A'),
 \qquad
 D'_k=B'Q_k.
\]
All these operators commute, and therefore $D'_k\geq0$.  Furthermore,
\[
 D'_k\uparrow B'P'_\Phi=B',
 \qquad
 0\leq D'_k\leq Q_k\leq kA'.
\]
The CP maps
\[
 \Psi_k(x)=V^*D'_k\pi(x)V
\]
therefore satisfy
$\Psi_k\uparrow\Psi$ and
$\Psi_k\leq_{\mathrm{cp}}k\Phi$.  Thus $\Psi\ll\Phi$.

Conversely, suppose that $\Psi\ll\Phi$.  Choose
$\Psi_k\uparrow\Psi$ and numbers $\alpha_k>0$ such that
$\Psi_k\leq_{\mathrm{cp}}\alpha_k\Phi$.  Since
$\Psi_k\leq_{\mathrm{cp}}\Psi\leq_{\mathrm{cp}}\Gamma$, let
$D'_k\in\pi(\cM)'_+$ be the Radon--Nikodym derivative of $\Psi_k$
relative to $\Gamma$.  The order isomorphism yields
\[
 0\leq D'_k\leq B',
 \qquad
 D'_k\leq\alpha_kA'.
\]
Moreover, $D'_k\uparrow B'$ strongly.  Since
$s(D'_k)\leq s(A')=P'_\Phi$ for every $k$, it follows that
$s(B')=P'_\Psi\leq P'_\Phi$.

Finally, recall that
\[
 (\Phi:\Psi)(x)=V^*(A':B')\pi(x)V.
\]
Here
\[
 A':B'=A'B',
\]
because $A'+B'=I$.  By injectivity of the Radon--Nikodym
correspondence,
\begin{align*}
 \Phi\perp\Psi
 &\iff \Phi:\Psi=0\\
 &\iff A'B'=0\\
 &\iff P'_\Phi P'_\Psi=0,
\end{align*}
which proves the singularity criterion.
\epf

\bex[Recovery of Ando's Lebesgue-type decomposition of bounded
positive operators]
\label{ex:ando-recovery}
Let $\cM=\mathbb C$, let $\cH=L^2(\cN)$ be a standard Hilbert
space for a von Neumann algebra $\cN$, and let $A,B\in\cN^+$.
Write
\[
 \Phi_A(z)=zA,
 \qquad
 \Phi_B(z)=zB,
 \qquad
 C=A+B,
\]
and put $p=s(C)$.

Under the canonical identification
$\mathbb C\odot\cH\simeq\cH$, the positive form associated with
$\Phi_C=\Phi_A+\Phi_B$ is
\[
 \gamma(\xi,\eta)=\ip{C\xi,\eta},
 \qquad \xi,\eta\in\cH.
\]
Let $\cH_\gamma$ be the Hilbert space obtained by separating and
completing $\cH$ with respect to $\gamma$, and denote by
$h_\gamma\colon\cH\to\cH_\gamma$
the canonical map.  Thus
\[
 \ip{h_\gamma(\xi),h_\gamma(\eta)}_{\cH_\gamma}
 =
 \ip{C\xi,\eta}_{\cH}.
\]

There is a canonical unitary
$U\colon\cH_\gamma\to p\cH$
determined by
\[
 Uh_\gamma(\xi)=C^{1/2}\xi,
 \qquad \xi\in\cH.
\]

A minimal Stinespring representation of $\Phi_C$ is given by
\[
 \pi(z)=zI_{\cH_\gamma},
 \qquad
 V\xi=h_\gamma(\xi)=U^*C^{1/2}\xi.
\]
Its Arveson Radon--Nikodym derivatives
$A',B'\in\bB(\cH_\gamma)_+$ are given by
\[
 V^*A'V=A,
 \qquad
 V^*B'V=B,
 \qquad
 A'+B'=I_{\cH_\gamma}.
\]
Recall that
\[
 (n\Phi_A:\Phi_B)(1)=nA:B.
\]
It follows that
\begin{align*}
 [\Phi_A]\Phi_B(1)
 &=
\lim_{n\to\infty}
 (n\Phi_A:\Phi_B)(1)
 \\
 &=
\lim_{n\to\infty}(nA:B)
 \\
 &=
 [A]B.
\end{align*}
Hence, we conclude that
\[
 [\Phi_A]\Phi_B=\Phi_{[A]B}.
\]
Thus the present construction recovers Ando's maximal
$A$-absolutely continuous part of $B$.

For completeness, assume first that $\ker A=\{0\}$ and consider
Kosaki's densely defined operator
\[
 T=B^{1/2}A^{-1/2},
 \qquad
 \operatorname{dom}T=\operatorname{ran}(A^{1/2}).
\]
Then
\[
 T^*=A^{-1/2}B^{1/2},
 \qquad
 \operatorname{dom}T^*
 =
 \{\xi\in\cH:
   B^{1/2}\xi\in\operatorname{ran}(A^{1/2})\}.
\]
Consequently, $T$ is closable if and only if
$\operatorname{dom}T^*$ is dense, which is equivalent to the
$A$-absolute continuity of $B$; see
\cite[Lemma~3]{kos_1984}.

We next consider the general case.  Since $A\leq C$, the map
\[
 C^{1/2}\xi\mapsto A^{1/2}\xi,
 \qquad \xi\in\cH,
\]
is well-defined and contractive, because
\[
 \|A^{1/2}\xi\|^2
 \leq
 \|C^{1/2}\xi\|^2.
\]
Since
$\operatorname{ran}(C^{1/2})$ is dense in $p\cH$, it extends
uniquely to a contraction
$H\colon p\cH\to\cH$
satisfying
\[
 A^{1/2}=HC^{1/2}.
\]

For every $\xi\in\cH$,
we have
\begin{align*}
 \ip{U^*H^*HUV\xi,V\xi}_{\cH_\gamma}
 &=
 \ip{H^*HC^{1/2}\xi,C^{1/2}\xi}_{\cH}
 \\
 &=
 \|HC^{1/2}\xi\|^2
 \\
 &=
 \|A^{1/2}\xi\|^2
 \\
 &=
 \ip{A\xi,\xi}_{\cH}
 \\
 &=
 \ip{A'V\xi,V\xi}_{\cH_\gamma}.
\end{align*}
Since $V\cH$ is dense in $\cH_\gamma$, 
we have
\[
 A'=U^*H^*HU.
\]
In particular,
\[
 \ker A'=U^*(\ker H).
\]

Since $A'+B'=I_{\cH_\gamma}$, the support criterion proved above
gives
\begin{align*}
 B\ll A
 &\iff
 \Phi_B\ll\Phi_A
 \\
 &\iff
 s(B')\leq s(A')
 \\
 &\iff
 \ker A'=\{0\}
 \\
 &\iff
 \ker H=\{0\}
 \\
 &\iff
 H\colon p\cH\to\cH
 \text{ is injective}.
\end{align*}
This is the support-reduced form of the injectivity criterion in
\cite[Corollary~13]{kos_2023}.

If $H$ is extended to an operator
\[
 \widehat H\coloneqq Hp\in\bB(\cH),
\]
then
the preceding condition is equivalently written as
\[
 B\ll A
 \iff
 \ker\widehat H=\ker C.
\]
When $\ker C=\{0\}$, we have $p=I_{\cH}$ and
$\widehat H=H$, so this reduces to the injectivity of $H$ on
all of $\cH$.

Finally,
\begin{align*}
 \Phi_A\perp\Phi_B
 &\iff
 [\Phi_A]\Phi_B=0
 \\
 &\iff
 \Phi_{[A]B}=0
 \\
 &\iff
 [A]B=0,
\end{align*}
which is precisely Ando's singularity condition.  Hence the
Stinespring support construction recovers both the absolutely
continuous and the singular parts of Ando's Lebesgue-type
decomposition.
\eex


\section*{Acknowledgements}

The author used Gemini and ChatGPT for English
language editing and mathematical information retrieval during the preparation of this manuscript.

\end{document}